\documentclass[12pt]{amsart}
\numberwithin{equation}{section}
\newtheorem{theorem}{Theorem}
\newtheorem{lemma}{Lemma}
\newtheorem{definition}{Definition}
\newtheorem{proposition}{Proposition}
\numberwithin{theorem}{section} \numberwithin{lemma}{section}
\newtheorem{corollary}{Corollary}
\numberwithin{proposition}{section}
\numberwithin{definition}{section}
\def\al{\aligned}
\def\eal{\endaligned}

\def\be{\begin{equation}}
\def\ee{\end{equation}}
\def\lab{\label}
\begin{document}

\tracingpages 1
\title[singular data]{\bf Local estimates on two linear
parabolic equations with singular coefficients}
\author{Qi S. Zhang}
\address{Department of Mathematics\\
University of California Riverside, Riverside, CA 92521}
\date{Pacific J. Math. Vol.223, No2, 2006, p367-396}

\begin{abstract}
We treat the heat equation with singular drift terms and its
generalization: the linearized Navier-Stokes system. In the first
case, we obtain boundedness of weak solutions for highly singular,
"supercritical" data. In the second case, we obtain regularity
result for weak solutions with mildly singular data ( those in the
Kato class). This not only extends some of the classical
regularity theory in [AS], [CrZ] and others from the case of
elliptic and heat equations to that of linearized Navier-Stokes
equations but also proves an unexpected gradient estimate, which
extends the recent interesting boundedness result [O].
\end{abstract}

\maketitle

\tableofcontents
\section{Introduction}

 The goal of the
paper is to prove local boundedness and other regularity of weak
solutions to the next two parabolic equations.

\[
\Delta u(x, t) - b(x, t) \nabla u(x, t) - \partial_t u(x, t)=0,
\quad (x, t) \in \Omega \subset {\bf R}^n \times {\bf R},
\leqno(1.1)
\]
\[
\begin{cases}
\Delta u(x, t) - b(x, t)  \nabla u(x, t) + \nabla P(x, t) -
\partial_t
u(x, t)=0, \quad (x, t) \in \Omega \subset {\bf R}^3 \times {\bf R},\\
div u = 0, \ div b=0, \ b(\cdot, t) \in L^2_{loc}.
\end{cases}
\leqno(1.2)
\]Here $\Delta$ is the standard Laplacian and $b=b(x, t)$ is a  given
$L^2_{loc}$ singular vector field  to be specified later. $\Omega$
is a domain.

There has been a mature theory of existence and regularity for
equation (1.1) (see [LSU], [Lieb] e.g.). For instance when
$b=b(x)$ and $|b| \in L^p_{loc}({\bf R}^n)$, $p>n$, weak solutions
to (1.1) are locally bounded and H\"older continuous. This
condition is sharp in general. Here is an example (see [HL] p108).
The function $u=\ln \ln |x|^{-1} - \ln \ln R^{-1}$ is an unbounded
weak solution of
\[
\Delta u + b \nabla u =0
\]
in the ball $B(0, R)$ in ${\bf R}^2$, $R<1$. Here $b = \nabla u =
- \frac{\nabla |x|}{|x| \ln |x|^{-1}}$ and hence $b \in L^2_{loc}$
with $n=2$.

The first goal of the paper is to show that the simple condition
$div b \le 0$ will ensure that weak solutions of (1.1) are locally
bounded when the data $b$ is  almost twice as singular as before.
This will be made precise in Theorem 1.1 and Remark 1.1 below.
Thus one has achieved a leap in boundedness condition rather than
a marginal improvement.

Clearly  a strong impetus still exists for the study of parabolic
equations with very singular coefficients. In the study of
nonlinear equations with gradient structure such as the
Navier-Stokes equations and harmonic maps, highly singular
functions occur naturally. So, it is very important investigate a
possible gain of regularity in the presence of singular drift term
$b$. This line of research has been followed in the papers [St],
[KS], [Os], [CrZ], [ChZ], [CS] and [Se]. Under the condition $|b|
\in L^n({\bf R}^n)$, Stampacchia [St] proved that bounded
solutions of $\Delta u + b \nabla u=0$ are H\"older continuous. In
the paper [CrZ], Cranston and Zhao proved that solutions to this
equation are continuous when $b$ is in a suitable Kato class i.e.
$\lim_{r \to 0} \sup_x \int_{|x-y| \le r}
\frac{|b(y)|}{|x-y|^{n-1}} dy =0$. In the paper [KS] Kovalenko and
Semenov proved the H\"older continuity of solutions to (1.1), when
$| b|^2$ is independent of time and  is sufficiently small in the
form sense, i.e., for a sufficiently small $\epsilon>0$,
\[
\int_{{\bf R}^n}  |b|^2(x)  \phi^2(x) dx \le \epsilon
\int_{\mathbf{R}^n} |\nabla \phi|^2(x) dx, \quad \phi\in
C^\infty_0(\mathbf{R}^n).
\]It is a well known fact that form boundedness condition provides
a more general class of singular functions than corresponding
$L^p$ class, Morrey-Campanato class and Kato class functions.
  This result was  generalized in [Se] to
equations with leading term in divergence form. In [Os], Osada
proved, among other things, that the fundamental solution of (1.1)
has global Gaussian upper and lower bound when $\bf b$ is the
derivative of bounded functions (in distribution sense) and $div
b=0$. Recently in the paper [LZ], H\"older continuity of solutions
to (1.1) was established when $b=b(x)$, $|b|^2$ is form bounded
and $div {b} =0$. Most recently, in [Z2], we considered (1.1) with
time-dependent functions $b=b(x, t)$.  It was proven that weak
solutions to (1.1) are locally bounded provided that $div b =0$
and for a fixed $m>1$, $|b|^m$ is form bounded. That is for any
$\phi \in C^\infty({\bf R}^n \times (0, \infty) )$ with compact
support in the spatial direction,
\[
\int \int_{{\bf R}^n} |b| ^m \phi^2 dxdt \le k \int \int_{{\bf
R}^n} |\nabla \phi|^2 dxdt
\]where $k$ is independent of $\phi$. Note the key improvement
over previous result is that the power on $b$ drops from $2$ to
any number greater than $1$.

It is interesting to note that this class of data $b$ contains the
velocity function in the $3$ dimensional Navier-Stokes equations.
As a result we gave a different proof of the {\it local}
boundedness of velocity in $2$ dimensional case. Moreover assuming
a {\it local}  bound in the pressure, we prove boundedness of
velocity in $3$ the dimensional case.

The first goal of the paper is to treat the end point case of the
above condition i.e. $m=1$. We will prove that weak solutions to
(1.1) are locally bounded provided that $|b| [\ln (1+ |b|)]^2$ is
form bounded and $div b \le 0$.

Many authors have also studied the regularity property of the
related heat equation $ \Delta u + V u - u_t =0$. Here $V$ is s
singular potential. We refer the reader to the papers by Aizenman
and Simon [AS], Simon [S] and the reference therein. The function
$V$ is allowed in the Kato class which is a little more singular
than the corresponding $L^p$ class. It remains a  challenging
problem is to push this theory to broader class of functions.

In this paper we use the following definition of weak
 solutions.

{\it {\bf Definition 1.1} Let $D \subseteq {\bf R}^n$ be a domain
and $T \in (0, \infty]$. A function $u$ such that $u, |\nabla u|
\in L^2_{loc}(D \times [0, T])$ is a weak solution to (1.1) if:
for any $\phi \in C^{\infty}_0(D \times (-T, T))$, there holds
\[
\int^T_0\int_D ( u \partial_t \phi - \nabla u \nabla \phi)
 dxdt - \int^T_0\int_D b \nabla u \  \phi \  dxdt
 = - \int_D u_0(x) \phi(x, 0) dx.
\]}

\begin{theorem}
\label{th:1.1} Suppose $div b \le 0$ in the weak sense, $b \in
L^2_{loc}$,  and that $|b| [\ln (1 + |b|)]^2$ is form bounded i.e.
there exists $k>0$ such that for any $\phi \in C^\infty({\bf R}^n
\times (0, \infty) )$ with compact support in the spatial
direction,
\[
\int \int_{{\bf R}^n} |b| [\ln (1 + |b|)]^2 \phi^2 dxdt \le k \int
\int_{{\bf R}^n} |\nabla \phi|^2 dxdt. \leqno(1.3)
\]
Then weak solutions to equation (1.1) are locally bounded.
\end{theorem}
\medskip

{\it Remark 1.1.} In the special case that $b$ is independent of
time and  $b \in L^p({\bf R}^n)$ with $p>n/2$, then it is easy to
check that (1.3) is satisfied. Recall that the standard theory
essentially only allows functions in $L^p$ with $p
> n$. The strength of the theorem comes from the fact that weak
solutions are locally bounded in any domain regardless of its
value on the parabolic boundary. If the domain is ${\bf R}^n
\times (0, \infty)$ or if the initial Dirichlet boundary condition
is imposed, then using a Nash type estimate, one can show that
solutions are locally bounded when $t>0$ as long as the
fundamental solution is well defined. In this case one can choose
$b$ to be as singular as any $L^2$ functions. Actually the
presence of $b$ is totally irrelevant except for the purpose of
making the integrals in the definition of a weak solution finite.
Here is a sketch of the proof. Let u solves
\[
\begin{cases}
 \Delta u - b \nabla u - u_t =0,  \qquad \text{in} \
 \text D \times (0, \infty),\\
u(x, t) = 0, \qquad (x, t) \in \partial D \times (0, \infty)\\
  u(x, 0) = u_0(x).
\end{cases}
  \]
Let $G(x, t; y, t)$ be the fundamental solution with initial
Dirichlet boundary condition. If $div b=0$, differentiating in
time shows that $\int_D G(x, t; y, 0) dy \le 1$. By Nash
inequality, one has
\[
\frac{d}{dt} \int_D G(x, t; y, t)^2 dy  = - 2 \int_D |\nabla G(x,
t; y, t)|^2 dy  \le - c \frac{ \big{(} \int_D G(x, t; y, t)^2 dy
\big{)}^{1+(2/n)}}{\big{(} \int_D G(x, t; y, 0) dy \big{)}^{4/n}}.
\]
Hence
\[
G(x, t; y, 0) \le c/t^{n/2}.
\]Therefore
\[
u(x, t) = \int_D G(x, t; y) u_0(y) dy
\] is bounded as soon as $t>0$
and $u_0$ is in $L^1(D)$.

 However this does not imply local
boundedness of weak solutions unconditionally. It would be
interesting to establish existence and more regularity result for
(1.1) with the singular data in Theorem 1.1.

We should mention that in the time independent elliptic case,
there was already strong indication that standard regularity
theory can be improved in the presence of divergence free data. In
the important papers [FR1-2], local boundedness of Green's
function (away from the singularity) of the operator $ -\Delta  -
b \nabla $ with Dirichlet boundary conditions was proved, under
the conditions:  $n=5$, $|\nabla^2 b| \in L^{4/3}$, $b=0$ on the
boundary and $div b =0$ ( c.f. Lemma 1.48 and Lemma1.11 [FR1]).
The upshot of the result is that the bounds on the Green's
function is independent of the norm $|\nabla^2 b| \in L^{4/3}$.
The proof uses essentially the fact that the Green's function
vanishes on the boundary. So the drift term is integrated out. In
contrast we do not have the benefit of a zero boundary.

In the three dimension case, we derive a further regularity
result:

\begin{corollary}
%col1.
 Assume $|b| \in L^{\infty}([0, T], L^2({\bf R}^3))
 \cap L^q(
 {\bf R}^3 \times [0, T])$ with $q>3$ and $div b =0$. Suppose $u$ is a weak solution of (1.1) in ${\bf R}^3
\times [0, T]$ with $\int_{{\bf R}^3} u^2(x, 0) dx < \infty$. Then
$u$ is locally bounded and for almost every $t$, $u(\cdot, t)$ are
H\"older continuous. \proof
\end{corollary}

According to Corollary 1 in [Z2], $b $ satisfies the conditions of
Theorem 1.1. For completeness we provide the proof here.

 Let us take $m=4/3$ and
$p=2/m=3/2$. Then, by H\"older's inequality,
\[
\aligned
 \int^T_0 \int_{{\bf R}^3}& |b |^{4/3} \phi^2 dxdt \le \int^T_0
\bigg{(} \int_{{\bf R}^3} |b |^{mp} dx \bigg{)}^{1/p}  \ \bigg{(}
\int_{{\bf R}^3}
\phi^{2p/(p-1)} dx \bigg{)}^{(p-1)/p} dt\\
&=\int^T_0 \bigg{(} \int_{{\bf R}^3} |b |^2 dx \bigg{)}^{2/3} \
\bigg{(} \int_{{\bf R}^3}
\phi^6 dx \bigg{)}^{1/3} dt\\
&\le \sup_{t \in [0, T]} \bigg{(} \int_{{\bf R}^3} |b |^2(x, t) dx
\bigg{)}^{2/3}
\ \int^T_0 \bigg{(} \int_{{\bf R}^3} \phi^6 dx \bigg{)}^{1/3} dt\\
&\le C \sup_{t \in [0, T]} \bigg{(} \int_{{\bf R}^3} |b |^2(x, t)
dx \bigg{)}^{2/3} \ \int^T_0 \int_{{\bf R}^3} |\nabla \phi|^2 dx
dt.
\endaligned
\]The last step is by Sobolev imbedding. This shows that condition
(1.3) holds.

Hence Theorem 1.1 implies that $u$ is locally bounded.  Notice the
fact that
\[
\int_{{\bf R}^3} u^2(x, t) dx
\]is non-increasing in time since $b$ is divergence free.
By this and Theorem 1.1, we know that $u \in L^{\infty}({\bf R}^3
\times [t_0, T])$ for any $t_0>0$.

Denote by $G_0$ the Gaussian heat kernel of the heat equation.
Then, for $t>t_0$,
\[
u(x, t) = \int_{{\bf R}^3} G_0(x, t; y, t_0) u(y, t_0) dy -
\int^t_{t_0} \int_{{\bf R}^3} G_0(x, t; y, s) b  \nabla u(y, s)
dyds.
\]Since $b$ is divergence free, we have
\[
u(x, t) = \int_{{\bf R}^3} G_0(x, t; y, t_0) u(y, t_0) dy +
\int^t_{t_0} \int_{{\bf R}^3} \nabla_y G_0(x, t; y, s) b u(y, s)
dyds.
\]Therefore, in the weak sense,
\[
\aligned
 \nabla_x u(x, t) &= \int_{{\bf R}^3} \nabla_x G_0(x, t; y,
t_0) u(y, t_0) dy + \int^t_{t_0} \int_{{\bf R}^3} \nabla_x
\nabla_y G_0(x, t; y, s) b u(y, s) dyds\\
& \equiv I_1(x, t) + I_2(x, t).
\endaligned
\]It is well known that
\[
\nabla_x \nabla_y G_0(x, t; y, t_0) =-\nabla_x \nabla_x G_0(x, t;
y, t_0)
\]is a parabolic Calderon-Zygmond kernel (see [Lie] e.g.).
Hence by our assumption on $b$ and the fact that $u$ is bounded in
${\bf R}^3 \times [t_0, T]$, the second term $I_2$  in the last
integral is in $L^q({\bf R}^3 \times [t_0, T])$, $q>3$. It follows
that
\[
|I_2(\cdot, t) | \in L^q({\bf R}^3), \qquad q>3.
\]for a.e. $t$. Sobolev imbedding theorem then shows that
 $u (\cdot, t)$ is
H\"older continuous for a.e. $t$. \qed

\medskip

 Next we turn to equation (1.2), which is the first step in
tackling the full Navier-Stokes equations. When $b=0$, (1.2) is
just the Stokes equations which has been studied for long time.
Our focus is on how to allow $b$ as singular as possible while
retaining the boundedness of weak solutions. As far as equation
(1.2) is concerned, our result does not improve the standard
theory as dramatically as for equation (1.1). We have to restrict
the data $b$ in a suitable Kato class for (1.2). Nevertheless,
Theorem 1.2 still generalizes the key part of the important work
[AS], [CrZ] on the elliptic equations to the case of linearized
Navier-Stokes system. Moreover, we even obtain gradient estimates
for solutions of (1.2) while only continuity was expected.

As pointed out in the papers [AS] and [Si], Kato class functions
are quite natural objects in studying elliptic and parabolic
equations with singular lower order terms. Roughly speaking, a
function is in a Kato class with respect to an equation if its
convolutions with certain kernel functions are small in some
sense. The kernel function usually is related to the fundamental
solution of the principal term of the equation. For instance, for
the equation $\Delta u(x) + V(x) u(x) =0$ in ${\bf R}^n$, $n \ge
3$, the function $V$ is in Kato class if
\[
\lim_{r \to 0} \sup_{x} \int_{B(x, r)} \frac{|V(y)|}{|x-y|^{n-2}}
dy = 0. \leqno(1.4)
\]In [AS], it is proven that weak solutions to $\Delta u + V u=0$
are continuous and satisfy a Harnack inequality when $V$ is in the
above Kato class. Numerous papers have been written on this
subject in the last thirty years, mainly in the context of
elliptic and heat equations.

In the context of Navier-Stokes equations, the corresponding time
dependent Kato class was defined recently in [Z3], which mirrors
those for the heat equation [Z]. Normally, with data in the Kato
class, weak solutions of elliptic equations are just continuous as
proven in [AS], [CrZ]. It was proved  in [Z3] that weak solutions
to (1.2) are bounded when $b$ is in the Kato class. Here we prove
that the {\it spatial gradient} of solutions to (1.2) are bounded
provided that $b$ is in the Kato class locally. Let us mention
that one can use the idea of Kato class to recover some (but not
all) the decay estimate in the interesting papers [Scho1-2] and to
prove some pointwise decay estimate (see [Z3]).

In order to make our statement precise, we introduce a number of
notations. Throughout the paper, we write
\[
K_1(x, t; y, s) = \begin{cases} \frac{1}{[ \ |x-y| + \sqrt{t-s} \
]^{n+1}}, \quad t \ge  s, x \neq y\\
0, s<t.
\end{cases}
\leqno (1.5)
\]We write $Q_r(x, t) = B(x, r) \times [t-r^2, t]$.

\medskip

{\it {\bf Definition 1.2}. A vector valued function $b=b(x, t) \in
L^1_{loc}({\bf R}^{n+1})$ is in class $K_1$ if it satisfies the
following  condition:
\[
\lim_{h \to 0 } \sup_{(x, t) \in {\bf R}^{n+1}} \int^{t}_{t-h}
\int_{{\bf R}^n} [K_1(x, t; y, s) + K_1(x, s; y, t-h) ] |b(y, s)|
dyds =0, \leqno (1.6)
\]}
\medskip

For clarity of presentation, given $t>l$,  we introduce the
quantity
\[
B(b, l, t) = \sup_{x} \int^{t}_l \int_{{\bf R}^n} [K_1(x, t; y, s) +
K_1(x, s; y, l) ] |b(y, s)| dyds. \leqno(1.6')
\]

By the example in Remark 1.2 in [Z3], we see that the function
class $K_1$ permits solutions which  are very singular. In case
the spatial dimension is $3$,  a function in this class can have
an apparent singularity of certain type that is not $L^p_{loc}$
for any $p>1$ and of dimension $1$. One can also construct time
dependent functions in $K_1$ with quite singular behaviors. The
class $K_1$ also contains the space $L^{p,q}$ with $n/p + 2/q<1$,
which sometimes is referred to as the Prodi-Serrin class. For the
nonlinear Navier-Stokes equation, if a weak solution is known to
be in this class, then it is actually smooth. As for the
linearized equation (1.2), following the argument in [Ser], it is
clear that weak solutions are bounded if $b$ is in the above
$L^{p, q}$ class. Now we are able to prove that the spatial
gradient of weak solutions are bounded {\it automatically},
without resorting to the nonlinear structure. Using H\"older's
inequality, one can see that the class $K_1$ also contains the
Morrey type space introduced in [O] by O' Leary, where boundedness
of weak solutions in Morrey space are proven.

One can also define a slightly bigger Kato class by requiring the
limit in (1.6) to be a small positive number rather than $0$. We
will not seek such generality this time. The appearance of two
kernel functions is due to the asymmetry of the equation in time
direction.

\medskip

Let $D$ for a domain in ${\bf R}^3$ and $T>0$. Following standard
practice, we will use this definition for solutions of (1.2)
throughout the paper.

\medskip

\noindent{\it  {\bf Definition 1.3.} A divergence free vector
field $u \in L^{\infty}(0, T; L^2(D)) \cap L^2(0, T; W^{1, 2}(D))$
is called a
 (weak) solution of (1.2) if:

for any vector valued $\phi \in C^{\infty}(D \times [0, T])$ with
$ div \phi =0$ and $\phi=0$ on $\partial D \times [0, T]$, $u$
satisfies
\[
\int^{t_2}_{t_1}\int_{{\bf R}^n} <u, \partial_t \phi + \Delta
\phi>
 dxdt - \int^{t_2}_{t_1}\int_{{\bf R}^n} <b \nabla u,  \phi>  dxdt
 = - \int_{{\bf R}^n} <u(x, t), \phi(x, t)> |^{t_2}_{t_1} dx.
\]}

Next we state the theorem on equation (1.2), the linearized
Navier-Stokes equations.

\bigskip

\begin{theorem}
\label{th:1.2}

Let $u$ be a solution of (1.2) in a domain $\Omega \subset {\bf
R}^3 \times {\bf R}$.  Suppose $Q_{4r}(x, t) \subset \Omega$, $div
b =0$ and that $b |_{Q_{2r}(x, t)}$ is in class $K_1$ and $b \in
L^2_{loc}$. Then both $u$ and $|\nabla u|$ are bounded functions
in $Q_{2r}(x, t)$.

Moreover, for some positive constants $C=C(b)$ and $r_0$,
depending on the size of the Kato norm of $b$,   there hold, when
$0<r<r_0$,
\[
\aligned
 |u(x, t)| &\le  \frac{C}{r^5} \int_{Q_{2r}(x, t)} |u(y,
s)| dyds,\\
|\nabla u(x, t)| &\le \frac{C}{r^5} \int_{Q_{2r}(x, t)} |\nabla
u(y, s)| dyds + \frac{C}{r^6} \int_{Q_{2r}(x, t)} |u(y,
s)-\overline{u}_{Q_{2r}}| dyds.
\endaligned
\leqno(1.7)
\]Here $\overline{u}_{Q_{2r}}$ is the average of $u$ in $Q_{2r}(x, t)$.

If in addition that $u(\cdot, t) \in W^{1, 2}_0(B(x, 2r))$ for
a.e. $t$, then
\[
|\nabla u(x, t)| \le \frac{C}{r^5} \int_{Q_{2r}(x, t)} |\nabla
u(y, s)| dyds. \leqno(1.8)
\]
\end{theorem}

\medskip

{\it Remark 1.2.} One may think that the last term on the right
hand side of the gradient estimate (1.7) is too singular to be
natural, especially when $r \to 0$. However, both (1.7) and
equation (1.2) are actually scaling invariant under the scaling:
$u_\lambda = u(\lambda x)$, $b_\lambda = \lambda b(\lambda x)$. So
(1.7) is in the right setting. Moreover, the gradient estimate in
(1.7) immediately simplifies to (1.8) when $u$ vanishes on the
lateral boundary or when $u$ enjoys a suitable extension property
in $W^{1, 2}$ space.

{\it Remark 1.3.} The reader may wonder whether any extra
regularity in the time direction is possible. The answer is no as
indicated in the example in [Ser], restated in [O]. Let $\phi$ be
a harmonic function in ${\bf R}^3$ and $a=a(t)$ be an integrable
function. Then $u= a \nabla \phi$ is a weak solution of the
Navier-Stokes equation. Obviously in the time direction $u$ is no
more regular than $a$.

\medskip

Theorem 1.1 and 1.2 will be proven in section 3 and 4
respectively. The idea for the proof of Theorem 1.2 is to combine
a recent localization argument in [O] with a refined iteration.
Using the idea in the proof of Theorem 1.1,  in Section 4  we
introduce a sufficient condition on the velocity that implies
boundedness of weak solutions of 3 dimensional Navier-Stokes
equations. The main improvement is that no absolute value of the
velocity is involved.

\section{Proof of Theorem 1.1}
%2

Since the drift term $b$ in (1.1) can be much more singular than
those allowed by the standard theory, the existence and uniqueness
of  weak solutions of (1.1) can not be taken for granted. In order
to proceed first we need  some approximation results whose proof
can be found in [Z2], with only small modifications.

\begin{proposition}
% Prop2.1
Suppose $u$ is a weak solution of  equation (1.1) in the cube $Q =
D \times [0, T]$, where $b$ satisfies the condition in Theorem
1.1, part (a). Here $D$ is a domain in ${\bf R}^n$.  Then $u$ is
the $L^1_{loc}$ limit of functions $\{ u_k \}$. Here $\{ u_k \}$
is a weak solution of (1.1) in which $b$ is replaced by smooth
$b_k$  such that $div b_k \le 0$ and  $b_k \to b$ strongly in
$L^2(Q)$, $k \to \infty$. \proof
\end{proposition}

The proof is almost identical to that of Proposition 2.4 in [Z2].
The only difference is that we are assuming $div b \le 0$ instead
of $div b =0$. Let $G=G(x,t; y, s)$ be the fundamental solution of
(1.1) with $div b \le 0$ and $b$ smooth. Then it is easy to show
by differentiation that
\[
\int_{{\bf R}^n} G(x, t; y, s) dx \le 1.
\]The rest of the proof is identical to that of Proposition 2.4
[Z2]. \qed

\medskip

Now we are ready to give a

{\bf Proof of Theorem 1.1.}

\medskip

{\bf Step 1.} $L^2$ gradient estimate. \medskip

By the above approximation result, we can and do assume that the
vector $b$ is smooth and we will differentiate freely as we wish.

In this section, we actually will prove the following local "mean
value" property.

Let $u$ be a  nonnegative solution of (1.1) in the parabolic cube
$Q_{\sigma_0 r} =B(x, \sigma_0 r) \times [t-(\sigma_0 r)^2, t]$.
Here $x \in {\bf R}^n$,
  $r>0$, $t>0$ and $\sigma_0$ is a suitable number greater than $1$.
  Suppose $b$ satisfies Condition A in
  $Q_{\sigma_0 r}$. Then there exist $C=C(r, b)>0$ such that
\[
\sup_{Q_r} u^2 \le C(r, b) \frac{1}{|Q_{\sigma_0 r}|} \int_{Q_{
\sigma_0 r}} u^2 dyds.
\]

Pick a solution $u$ of (1.1) in the parabolic cube $Q_{\sigma r}
=B(x, \sigma r) \times [t-(\sigma r)^2, t]$, where $x \in {\bf
R}^n$, $\sigma>1$,  $r>0$ and $t>0$. By direct computation, for
any rational number $p \ge 1$, which can be written as quotient of
two integers with the denominator being odd, one has
\[
\Delta u^p - b \nabla u^p  - \partial_t u^p = p(p-1) |\nabla u|^2
u^{p-2}. \leqno (2.1)
\]Here the condition on $p$ is to ensure that $u^p$ makes sense
when $u$ changes sign. Actually $u^2$ is a sub-solution to (1.1).
Hence one can also assume that $u$ is a non-negative sub-solution to
(1.1) by working with $u^2$.

Choose $\psi=\phi(y) \eta(s)$ to be a refined cut-off function
satisfying
\[
supp \  \eta \subset [t- (\sigma r)^2, t]; \quad \eta(s)=1, \quad
s \in [t- r^2, t]; \quad |\eta'| \le 2/( (\sigma-1) r)^2; \quad  0
\le \eta \le 1;
\]
\[
supp \  \phi \subset B(x, \sigma r); \quad \phi(y)=1, \quad y \in
B(x, r);  \quad 0 \le \phi \le 1;
\]
\[
\frac{|\nabla \phi|}{\phi} \le  \frac{A}{(\sigma-1) r} |\ln \phi
|^{3/2}, \qquad A>0.
\]
By modifying the following function
\[
\exp \big{(} -\frac{ \sigma^2}{\sigma^2 - |x-y|^2} \big{)}^k
\]
and scaling, it is easy to show that such a function $\phi$
exists. Here $k$ is a sufficiently large number.

Denoting $w=u^p$ and using $w \psi^2$ as a test function on (2.1),
one obtains
\[
\int_{Q_{\sigma r}} (\Delta w - b \nabla w  -
\partial_s w) w \psi^2 dyds  = p(p-1) \int_{Q_{\sigma r}}
|\nabla u|^2 w^2 u^{-2} \ge 0.
\]Using integration by parts, one deduces
\[
\int_{Q_{\sigma r}} \nabla(w \psi^2) \nabla w dyds \le -
\int_{Q_{\sigma r}} b \nabla w (w \psi^2) dyds  - \int_{Q_{\sigma
r}} (\partial_s w) w \psi^2 dyds. \leqno (2.2)
\]By direct calculation,
\[
\aligned &\int_{Q_{\sigma r}} \nabla(w \psi^2)  \nabla w dyds =
\int_{Q_{\sigma r}} \nabla [(w \psi )  \psi] \nabla w dyds\\
&=\int_{Q_{\sigma r}} [ \ \nabla (w \psi ) ( \ \nabla (w \psi ) -
(\nabla \psi) w ) + w \psi \nabla \psi \nabla w ] dyds\\
&=\int_{Q_{\sigma r}} [\  |\nabla (w \psi )|^2 -|\nabla \psi |^2
w^2 \ ]dyds.
\endaligned
\]Substituting this to (2.2), we obtain
\[
\aligned
\int_{Q_{\sigma r}} &|\nabla (w \psi )|^2 dyds \\
&\le - \int_{Q_{\sigma r}} b \nabla w (w \psi^2) dyds  -
\int_{Q_{\sigma r}} (\partial_s w) w \psi^2 dyds + \int_{Q_{\sigma
r}} |\nabla \psi |^2 w^2 dyds.
\endaligned \leqno (2.3)
\]

Next, notice that
\[
\aligned \int_{Q_{\sigma r}}& (\partial_s w) w \psi^2 dyds =
\frac{1}{2}
\int_{Q_{\sigma r}} (\partial_s w^2) \psi^2 dyds\\
&=-\int_{Q_{\sigma r}} w^2 \phi^2 \eta \partial_s \eta dyds +
\frac{1}{2} \int_{B(x, \sigma r)} w^2(y, t) \phi^2(y) dy.
\endaligned
\]Combining this with (2.3), we see that
\[
\aligned & \int_{Q_{\sigma r}} |\nabla (w \psi )|^2 dyds +
\frac{1}{2} \int_{B(x, \sigma r)} w^2(y, t) \phi^2(y) dy \\
&\le \int_{Q_{\sigma r}} (|\nabla \psi |^2 + \eta \partial_s \eta)
\ w^2 dyds -
\int_{Q_{\sigma r}} b (\nabla w) (w \psi^2) dyds  \\
&\equiv T_1+T_2.
\endaligned
\leqno(2.4)
\]The first term on the righthand side of (2.4) is
already in good shape. So let us estimate the second term as
follows.
\[
\aligned
T_2 &= - \int_{Q_{\sigma r}} b (\nabla w) (w \psi^2) dyds\\
&= - \frac{1}{2} \int_{Q_{\sigma r}} b \psi^2 \nabla w^2 dyds =
\frac{1}{2} \int_{Q_{\sigma r}} div (b \psi^2)  w^2 dyds\\
&=\frac{1}{2} \int_{Q_{\sigma r}} div b (\psi  w) ^2 dyds +
\frac{1}{2}
\int_{Q_{\sigma r}} b \nabla (\psi^2) w^2 dyds\\
&=\frac{1}{2} \int_{Q_{\sigma r}} div b (\psi  w) ^2 dyds +
\int_{Q_{\sigma r}} b (\nabla \psi) \psi w^2 dyds\\
&\le   \int_{Q_{\sigma r}} b (\nabla \psi) \psi w^2 dyds.
\endaligned
\]Here we just used the assumption that $div b \le 0$.

The next paragraph contains a key argument of the paper.

Let  $D>0$ be a number to be chosen later
\[
\aligned T_2 &\le  \int_{Q_{\sigma r}} |b| \  |\nabla \psi| \psi
w^2
dyds  \\
&\le \int_{|b| \ge D} |b| \  |\nabla \psi| \psi w^2 dyds
+\int_{|b| \le D} |b| \  |\nabla \psi| \psi w^2 dyds
\\
&\le \int_{|b| \ge D} |b| \  |\nabla \psi| \psi w^2 dyds + \frac{C
D }{(\sigma - 1) r} \int_{Q_{\sigma r}}  w^2 dyds.
\endaligned
\]Using
the property that $\psi = \phi \eta$ and
\[
|\nabla \phi| \le  \frac{A}{(\sigma -1) r} \phi | \ln \phi
|^{3/2},
\]
we have
\[
\aligned T_2& \le \int \int_{|b| \ge 1/\phi, |b| \ge D} |b|
|\nabla \phi| \phi w^2  dy \eta^2 ds  + \int \int_{|b|
\le 1/\phi}  |b| |\nabla \phi| \phi w^2 dy \eta^2 ds\\
&\qquad \qquad \qquad+
\frac{C D }{(\sigma - 1) r} \int_{Q_{\sigma r}}  w^2 dyds\\
&\le \frac{A}{(\sigma -1) r} \int \int_{|b| \ge 1/\phi, |b| \ge D}
|b| \ |\ln \phi|^{3/2} (\phi w)^2 dy \eta^2 ds \\
&\qquad \qquad + \frac{A}{(\sigma -1) r} \int \int_{B(x, \sigma
r), |b| \le 1/\phi} |b| \ |\ln \phi|^{3/2} (\phi w)^2  dy \eta^2
ds
\\
&\qquad \qquad \qquad+
\frac{C D }{(\sigma - 1) r} \int_{Q_{\sigma r}}  w^2 dyds\\
&\le \frac{A}{(\sigma -1) r} \int_{|b| \ge D} |b| (\ln |b|)^{3/2}
(\psi w)^2 dyds  + \frac{A}{(\sigma -1) r} \int_{Q_{\sigma r}}
\frac{1}{\phi} |\ln \phi|^{3/2} (\phi
w)^2 \eta^2 dyds\\
&\qquad \qquad \qquad +\frac{C D }{(\sigma - 1) r} \int_{Q_{\sigma
r}}  w^2 dyds\\
&\le \frac{A}{(\sigma -1) r (\ln D)^{1/2}} \int_{Q_{\sigma r}} |b|
(\ln |b|)^2 (\psi w)^2 dyds  + \frac{A}{(\sigma -1) r}
\int_{Q_{\sigma r}} \phi |\ln \phi|^{3/2}
w^2 \eta^2 dyds\\
&\qquad \qquad \qquad +\frac{C D }{(\sigma - 1) r} \int_{Q_{\sigma
r}}  w^2 dyds.\\
\endaligned
\]

By our assumptions on $b$,
\[
\int_{Q_{\sigma r}} |b| (\ln |b|)^2 (\psi  w)^2 dyds \le k
\int_{Q_{\sigma r}} |\nabla (\psi  w)|^2 dyds,
\]and the fact that $\phi |\ln \phi|^{3/2}$ is a bounded function, we
deduce
\[
T_2 \le \frac{A}{(\sigma -1) r (\ln D)^{1/2}} \int_{Q_{\sigma r}}
|\nabla (\psi w)|^2 dyds + \frac{C_1 D }{(\sigma - 1) r}
\int_{Q_{\sigma r}} w^2 dyds.
\]Now we choose $D$ so that $\frac{A}{(\sigma -1) r (\ln D)^{1/2}}
= \frac{1}{2}$, i.e. $D= e^{( 2A/ [(\sigma-1) r])^2}$, then
\[
T_2 \le  \frac{1}{2} \int_{Q_{\sigma r}} |\nabla (\psi  w)|^2 dyds
+ c_0 e^{( c_1/ [(\sigma-1) r])^2} \int_{Q_{\sigma r}}  w^2 dyds.
\leqno(2.6)
\]Here $c_0$ and $c_1$ are positive constants independent of $r$
and $\sigma$.

Combining  (2.4) with (2.6), we reach
\[
\int_{Q_{\sigma r}} |\nabla (w \psi )|^2 dyds +  \int_{B(x, \sigma
r)} w^2(y, t) \phi^2(y) dy  \le c_0 e^{( c_1/ [(\sigma-1) r])^2}
\int_{Q_{\sigma r}}  w^2 dyds. \leqno(2.7)
\]

\medskip

{\bf step 2.} $L^2 - L^\infty$ bounds.
\medskip

By modifying Moser's iteration moderately,  we deduce from (2.7)
the following $L^2-L^\infty$ estimate.
\[
\sup_{Q_r} u^2 \le C(r, b) \frac{1}{|Q_{\sigma_0 r}|}
\int_{Q_{\sigma_0 r}} u^2 dyds. \leqno(2.8)
\]Indeed, by H\"older's inequality,
\[
\aligned
 \int_{{\bf R}^n} &{(\phi w)}^{2(1+(2/n))} =
 \int_{{\bf R}^n} (\phi w)^2 \ (\phi w)^{4/n}\\
 &\le \big{(}
\int_{{\bf R}^n} {(\phi w)}^{2n/(n-2)} \big{)}^{(n-2)/n}
             \big{(}\int_{{\bf R}^n} {(\phi w)}^2 \big{)}^{2/n}.
\endaligned
\]Using  Sobolev inequality, one obtains
\[
\int_{{\bf R}^n} {(\phi w)}^{2(1+(2/n))} \le
            C \big{(}\int {(\phi w)}^2 \big{)}^{2/n}
             \big{(}\int_{{\bf R}^n}
            |\nabla (\phi w)|^2 \big{)}.
\]The last inequality, together with (2.7) implies, for some
$C_1>0$,
\[
\int_{Q_{\sigma' r }(x, t)}  u^{2p\theta} \le  \big{(} c_0 e^{(
c_1/ [(\sigma-\sigma') r])^2} \int_{Q_{\sigma r}(x, t)}
u^{2p}\big{)}^{\theta},
\]where $\theta = 1+ (2/n)$ and $\sigma'<\sigma$.

Take a  number $\rho>1$ so that $\rho^2< \theta$.
  We  set
$\tau_i=\rho^{-i}$, $\sigma_0=1/(1-\rho^{-1})$,
$\sigma_i=\sigma_{i-1}-\tau_i= \sigma_0-\Sigma^{i}_1 \tau_j$,
$p=\theta^i$, $i=1, 2, ...$. The above then yields, for some $c_2,
c_3>0$,
\[
\int_{Q_{\sigma_{i+1} r}(x, t)}  u^{2 \theta^{i+1} } \le c_3
\big{(} c^{i+1}_2 e^{ c^2_1 r^{-2} \rho^{2 i} } \int_{Q_{\sigma_i
r }(x, t)} u^{2 \theta^i}\big{)}^{\theta}.
\]After iterations the above implies, for some $c_4 >0$,
\[
\big{(}\int_{Q_{\sigma_{i+1} r }(x, t)}
        u^{2 \theta^{i+1}} \big{)}^{\theta^{-i-1}}
\le \exp (c_4 \Sigma^i_1 j \theta^{-j}) \ \exp (c_1 r^{-2}
\Sigma^i_1 \rho^{2j} \theta^{-j})
    \int_{Q_{\sigma_0 r}(x, t)} u^2.
\]Observe that $\rho^2/\theta<1$. Letting $i \rightarrow \infty$ and
observing that $\sigma_i \to 1$ as $i \to \infty$, we obtain
\[
\sup_{Q_r} u^2 \le C(r, b) \int_{Q_{\sigma_0 r}} u^2.
\]
 This completes
proof the theorem. \qed

\medskip

\section{Proof of Theorem 1.2}
%3

First we will need a short lemma concerning kernel function $K_1$
defined in (1.5). It was proved, among others things, in [Z3]. We
give a proof for completeness.

\begin{lemma}
%L3.1.
The following inequality holds for all $x, y, z \in {\bf R}^n$ and
$t> \tau >0$.
\[
\aligned K_1*b K_1 &\equiv \int^t_0\int_{\bf R^n} \frac{1}{(
|x-z|+\sqrt{t-\tau} )^{n+1}} \frac{|b(z, \tau)|}{(
|z-y|+\sqrt{\tau} )^{n+1}} dzd\tau \\
&\le C B(b, 0, t) K_1(x, t; y, 0). \endaligned \leqno(3.1)
\]Here, recalling from (1.6'),
\[
B(b, 0, t) \equiv \sup_{x \in {\bf R}^n} \int^{t}_0 \int_{{\bf
R}^n} [K_1(x, t; y, s) + K_1(x, s; y, 0) ] |b(y, s)| dyds. \leqno
(3.2)
\]
 \proof
\end{lemma}

Since
\[
|x-z| + \sqrt{t-\tau} + |y-z| + \sqrt{\tau} \ge |x-y| + \sqrt{t},
\]we have, either
\[
|x-z| + \sqrt{t-\tau} \ge \frac{1}{2} (|x-y| + \sqrt{t}),
\leqno(3.3)
\]or
\[
|z-y| + \sqrt{\tau} \ge \frac{1}{2} (|x-y| + \sqrt{t}).
\leqno(3.4)
\]

Suppose (3.3) holds then
\[
K_1*b K_1 \le \frac{ 2^n}{(|x-y| + \sqrt{t})^{n+1}}
\int^t_0\int_{\bf R^n} \frac{|b(z, \tau)|}{(|z-y|+\sqrt{\tau}
)^{n+1}} dzd\tau.
\]That is
\[
K_0*b K_1 \le \frac{ 2^n  B(b, 0, t)}{(|x-y| + \sqrt{t})^{n+1}}.
\leqno (3.5)
\]

Suppose (3.4) holds but (3.3) fails, then
\[
|z-y| + \sqrt{\tau} \ge \frac{1}{2} (|x-y| + \sqrt{t}) \ge |x-z| +
\sqrt{t-\tau}.
\]This shows
\[
\frac{1}{( |x-z|+\sqrt{t-\tau} )^{n+1} \ ( |z-y|+\sqrt{\tau}
)^{n+1}} \le \frac{2^n}{( |x-z|+\sqrt{t-\tau} )^{n+1} (
|x-y|+\sqrt{t} )^{n+1}}.
\]Substituting this to (3.1), we obtain
\[
K_1*b K_1 \le \frac{ 2^n}{(|x-y| + \sqrt{t})^{n+1}}
\int^t_0\int_{\bf R^n} \frac{|b(z, \tau)|}{(|x-z|+\sqrt{t-\tau}
)^{n+1}} dzd\tau.
\]That is
\[
K_1*b K_1 \le \frac{ 2^n  B(b, 0, t)}{(|x-y| + \sqrt{t})^n}.
\]

Clearly the only remaining case to consider is when both (3.3) and
(3.4) holds. However this case is already covered by (3.5). Thus
(3.1) is proven.  \qed

\medskip

Next we state and prove a representation formula for solutions of
(1.2) and their spatial gradient, following and extending the idea
in [O]. The formula for solutions is contained in [O]. However, we
will outline the proof since it is useful in the proof of the
formula for the gradient, which is a new contribution of this
paper.

\medskip

{\bf Remark 3.1.}  Let us note that at this moment, the
representation formula (3.6') below for the gradient is understood
as a comparison of two $L^1_{loc}$ functions in space-time. This
is legal for two reasons. First we assumed that $\nabla u$ is a
$L^2$ function a priori. Second, it is easy to check
 $K_1(\cdot,
\cdot; y, s)$ is $L^1_{loc}$ and $b \nabla u$ is $L^1_{loc}$ by
the assumption that $b$ is $L^2$. Therefore the last function on
the righthand side of (3.6') is a $L^1_{loc}$ function.

\medskip

\begin{lemma} (mean value inequality)
%L3.2.

(a). Let $u$ be a solution of (1.2) in the region $\Omega$.
Suppose $Q_{2r}(x, t) \subset \Omega$. Then there exists a
constant $\lambda$ such that
\[
\aligned
 |u(x, t)| &\le \lambda \frac{1}{r^5} \int_{Q_r(x,
t)-Q_{r/2}(x, t)} |u(y, s)| dyds \\
&\qquad + \lambda \int_{Q_r(x, t)} K_1(x, t; y, s) |b(y, s)| \
|u(y, s)| dyds. \endaligned
 \leqno(3.6)
\]

(b). Under the same assumption as (a), there exists a constant
$\lambda$ such that
\[
\aligned |\nabla u(x, t)| &\le  \frac{\lambda}{r^5} \int_{Q_r(x,
t)-Q_{r/2}(x, t)} | \nabla u(y, s)|  dyds\\
&\qquad  + \frac{\lambda}{r^6}
\int_{Q_r(x, t)-Q_{r/2}(x, t)} |u(y, s)| dyds\\
&\qquad  + \lambda \int_{Q_r(x, t)} K_1(x, t; y, s) |b(y, s)| \
|\nabla u(y, s)| dyds. \endaligned \leqno(3.6')
\]

\proof of (a).
\end{lemma}

Let $E=E(x, t; y, s)$ be the fundamental solution (matrix) of the
Stokes system in ${\bf R}^3 \times (0, \infty)$ and $E_k$ be the
$k$-th column of $E$. This function has been studied for a long
time. All of its basic properties we are using below can be found
in [So] and [FJR].  Fixing $(x, t)$, we construct a standard
cut-off function $\eta$ such that $\eta(y, s)=1$ in $Q_{r/2}(x,
t)$, $\eta(y, s)=0$ outside of $Q_r(x, t)$, $0 \le \eta \le 1$ and
$|\nabla \eta|^2 + |\Delta \eta| + |\partial_s \eta| \le c/r^2$.

Define a vector valued function
\[
\Phi_k=\Phi^{(x, t)}_k(y, s) = \frac{1}{4 \pi} curl \big{(} \eta
(y, s) \int_{{\bf R}^3} \frac{ curl E_k(x, t; z, s)}{|z-y|} dz
\big{)}. \leqno(3.7)
\]It is clear that when $t>s$,  $\Phi_k$ is a valid test function
for equation (1.2) since $\Phi_k$ is smooth, compactly supported
and divergence free. Using $\Phi_k$ as a test function on (1.2),
by Definition (1.2), we obtain
\[
\int^t_0\int u(y, s) \big{(} \Delta \Phi_k + b \nabla \Phi_k +
\partial_s \Phi_k \big{)} dyds = \lim_{s \to t} \int u(y, s)
\Phi_k(y, s) dy.
\]Here and later,we will suppress the superscript $(x, t)$ on
$\Phi$, unless there is a confusion.

Since $E_k$ is divergence free $curl \ curl E_k = - \Delta E_k$.
Thus
\[
\aligned \Phi_k(y, s) &= \eta(y, s) E_k(x, t; y, s) + \frac{1}{4
\pi} \nabla \eta(y, s) \times \int_{{\bf R}^3} \frac{ curl E_k(x,
t; z, s)}{|z-y|} dz\\
&\equiv \eta E_k + \overrightarrow{Z}.
\endaligned
\leqno(3.8)
\]
Using  the property of the fundamental matrix $E$ and the fact
that $\overrightarrow{Z}$ is a lower order term, it is easy to see
that
\[
\lim_{s \to t} \int u(y, s) \Phi_k(y, s) dy = \eta (y, t) u_k(x,
t) = u_k(x, t)
\]where $u_k$ is the k-th component of $u$. Hence
\[
\aligned u_k(x, t) &= \int_{Q_r(x, t)} u [ \Delta (\eta E_k) +
\partial_s (\eta
E_k)]dyds  + \int_{Q_r(x, t)} u [\Delta \overrightarrow{Z} +
\partial_s \overrightarrow{Z}] dyds\\
&\qquad \qquad +   \int_{Q_r(x, t)} u b \nabla \Phi_k dyds.
\endaligned
\]Here and later $u b \nabla \Phi_k = u \cdot \sum^3_{i=1} b_i
\partial_i \Phi_k$.

Note that $\Delta E_k + \partial_s E_k=0$ when $s<t$. The above
implies
\[
\aligned u_k&(x, t) \\
&= \int_{Q_r(x, t)} u(y, s) [ E_k(x, t; y, s) (\Delta \eta +
\partial_s \eta)(y, s) +
2 \nabla \eta(y, s) \nabla_y E_k(x, t; y, s)]dyds \\
&\qquad + \int_{Q_r(x, t)} u(y, s) [\Delta \overrightarrow{Z} +
\partial_s \overrightarrow{Z}](y, s) dyds\\
&\qquad \qquad+ \int_{Q_r(x, t)} u b \nabla_y \Phi_k(y, s) dyds\\
&\equiv J_1+J_2+J_3.
\endaligned
\leqno(3.9)
\]

We are going to estimate $J_1, J_2, J_3$ separately. By well known
estimates on $E$ (see [So] or [FJR] for example), for $(y, s) \in
Q_r(x, t)- Q_{r/2}(x, t)$,
\[
|E_k(x, t; y, s)| \le \frac{c}{(|x-y| + \sqrt{t-s})^3} \le
\frac{c}{r^3},
\]
\[
|\nabla_y E_k(x, t; y, s)| \le \frac{c}{(|x-y| + \sqrt{t-s})^4}
\le \frac{c}{r^4}.
\]Using these and the property of $\eta$, we see that
\[
|J_1| \le \frac{C}{r^5} \int_{Q_r(x, t)-Q_{r/2}(x, t)} |u(y, s)|
dyds. \leqno(3.10).
\]

Next we give an estimate $J_2$, which follows verbatim from $[O]$
p626. Recall that
\[
E_k(x, t; y, s) =  G(x, t; y, s) e_k + \frac{1}{4 \pi} \nabla
\partial_k \int_{{\bf R}^3} \frac{G(x, t; z, s)}{|y-z|} dz.
\]Here $G$ is the fundamental solution of the heat equation and
$\{ e_1, e_2, e_3 \}$ is the standard orthonormal basis of ${\bf
R}^3$. Since $curl \nabla =0$, from (3.8), $\overrightarrow{Z}$
defined in (3.8) takes a very simple form
\[
\overrightarrow{Z}= \frac{1}{ 4 \pi} \nabla_y \eta(y, s) \times
\bigg{(} \nabla_y \big{[}\int_{{\bf R}^3} \frac{ G(x, t; z,
s)}{|z-y|} dz \big{]} \times e_k \bigg{)}.
\]Hence
\[
J_2 = \frac{1}{ 4 \pi} \int_{Q_r(x, t)} u(y, s) [\Delta  +
\partial_s ] \bigg{[} \nabla_y \eta(y, s) \times
\bigg{(} \nabla_y \big{[}\int_{{\bf R}^3} \frac{ G(x, t; z,
s)}{|z-y|} dz \big{]} \times e_k \bigg{)} \bigg{]} dyds.
\leqno(3.11)
\]From here direct computation, using the estimates ([So], Chapter
2, Section 5),
\[
|D^l_s D^m_y \int_{{\bf R}^3} \frac{ G(x, t; z, s)}{|z-y|} dz| \le
\frac{C_{m, l}}{(|x-y| + \sqrt{t-s})^{1+m+2l}}, \leqno(3.12)
\] and the fact that $|x-y| \ge r/2$ here, it is easy to show that
\[
|J_2| \le  \frac{C}{r^5} \int_{Q_r(x, t)-Q_{r/2}(x, t)} |u(y, s)|
dyds. \leqno(3.13)
\]
Finally direct  computation using (3.12) shows that
\[
|J_3| \le \lambda \int_{Q_r(x, t)} K_1(x, t; y, s) |b(y, s)| \
|u(y, s)| dyds. \leqno(3.14)
\]Substituting (3.10), (3.13) and (3.14) to (3.9), we finish the
proof of part (a) of the lemma.
\medskip

{\bf Proof of part (b).}

Our next task is to prove the representation formula for $\nabla
u$. In (3.9) the cut-off function apparently depends on $(x, t)$.
In order to prove the gradient estimate, we need modify it a
little. For $(w, l) \in Q_{r/4}(x, t)$, we take
\[
\Phi_k=\Phi^{(w, l)}_k(y, s) = \frac{1}{4 \pi} curl \big{(} \eta
(y, s) \int_{{\bf R}^3} \frac{ curl E_k(w, l; z, s)}{|z-y|} dz
\big{)}.
\]as a text function for (1.2). Since $\eta(w, l)=1$ in this case,
following the computation before (3.9), we obtain
\[
\aligned u_k&(w, l) \\
&= \int_{Q_r(x, t)} u(y, s) [ E_k(w, l; y, s) (\Delta \eta +
\partial_s \eta)(y, s) +
2 \nabla \eta(y, s) \nabla_y E_k(w, l; y, s)]dyds \\
&\qquad + \int_{Q_r(x, t)} u(y, s) [\Delta \overrightarrow{Z} +
\partial_s \overrightarrow{Z}](y, s) dyds\\
&\qquad \qquad + \int_{Q_r(x, t)} u b \nabla_y \Phi^{(w, l)}_k(y,
s) dyds.
\endaligned
\leqno(3.9')
\]Here $\overrightarrow{Z}=\overrightarrow{Z}(w, l; y, s)$ is defined
 as $\overrightarrow{Z}$ in (3.8) except that $E_k(x, t; y, s)$ is
replaced by $E_k(w, l; y, s)$. i.e.
\[
\overrightarrow{Z} =\frac{1}{4 \pi} \nabla \eta(y, s) \times
\int_{{\bf R}^3} \frac{ curl E_k(w, l; z, s)}{|z-y|} dz.
\]We would like to differentiate (3.9) in the spatial variables.
However, since $\nabla u$ is only known as a $L^2$ function, we
have to consider the weak derivatives.

Let $\rho = \rho(w)$ be a smooth cut-off function supported in
$B(x, r/4)$. Then (3.9') implies, for $i=1, 2, 3$, and a.e. $l$,
\[
\aligned
&\int_{{\bf R}^3} u_k(w, l) \partial_{w_i} \rho(w) dw  \\
&= - \int_{{\bf R}^3} {\bf \bigg{(}}  \int_{Q_r(x, t)} u(y, s)
\partial_{w_i}[ E_k(w, l; y, s)
(\Delta \eta +
\partial_s \eta)(y, s) +
2 \nabla \eta(y, s) \nabla_y E_k(w, l; y, s)]dyds \\
&\qquad + \int_{Q_r(x, t)} u(y, s) \partial_{w_i}[\Delta
\overrightarrow{Z} +
\partial_s \overrightarrow{Z}](y, s) dyds\\
&\qquad \qquad +  \int_{Q_r(x, t)} u b \nabla_y
\partial_{w_i} \Phi^{(w, l)}_k(y, s) dyds {\bf \bigg{ )}} \rho(w) dw\\
&\equiv - \int_{{\bf R}^3} \bigg{(} M_1 + M_2 + M_3 \bigg{)}
\rho(w) dw.
\endaligned
\leqno(3.15)
\]Here we just used integration by parts which is legitimate since
we will show that $M_1, M_2$ are bounded functions and $M_3$ is
$L^1_{loc}$. Note also we should have integrated in the time
direction as well since $\nabla u$ is only known to $L^2$ in space
time. However, since all the estimates below are uniform for $l
\in [t-(r/4)^2, t]$, our estimates are valid.

Let us estimate $M_1$ first. Noting that $\partial_{w_i} E_k(w, l;
y, s) = -\partial_{y_i} E_k(w, l; y, s)$, we deduce
\[
M_1 = - \int_{Q_r(x, t)} u(y, s) [ \partial_{y_i} E_k(w, l; y, s)
(\Delta \eta +
\partial_s \eta)(y, s) +
2 \nabla \eta(y, s) \nabla_y \partial_{y_i} E_k(w, l; y, s)] dyds.
\]Using integration by parts, we obtain
\[
\aligned &M_1 =\\
&  \int_{Q_r(x, t)} \partial_{y_i} u(y, s) [ E_k(w, l; y, s)
(\Delta \eta +
\partial_s \eta)(y, s) +
2 \nabla \eta(y, s) \nabla_y  E_k(w, l; y, s)] dyds\\
&+ \int_{Q_r(x, t)} u(y, s) [  E_k(w, l; y, s)
\partial_{y_i} (\Delta \eta +
\partial_s \eta)(y, s) +
2 (\partial_{y_i} \nabla \eta(y, s)) \nabla_y  E_k(w, l; y, s)]
dyds.
\endaligned
\]By standard properties of  $E$ and the bounds on $\eta$ and
its derivatives, we have
\[
|M_1| \le \frac{c}{r^5} \int_{Q_r(x, t)-Q_{r/2}(x, t)} |\nabla
u(y, s)| dyds +\frac{c}{r^6} \int_{Q_r(x, t)-Q_{r/2}(x, t)} |u(y,
s)| dyds. \leqno(3.16)
\]Here we also utilize the fact that the arguments $(w, l)$ and $(y, s)$
have a parabolic distance of at least $r/4$

Next we need to find an upper bound for $M_2$. Recall from the
formula before (3.11) that
\[
\aligned \overrightarrow{Z}&= \frac{1}{ 4 \pi} \nabla_y \eta(y, s)
\times \bigg{(} \nabla_y \big{[}\int_{{\bf R}^3} \frac{ G(w, l; z,
s)}{|z-y|} dz
\big{]} \times e_k \bigg{)}\\
&=\frac{1}{ 4 \pi} \nabla_y \eta(y, s) \times \bigg{(} \nabla_y
\big{[}\int_{{\bf R}^3} \frac{ G(w, l; z+y, s)}{|z|} dz \big{]}
\times e_k \bigg{)}.
\endaligned
\]Using the vector identity $F \times (G \times H) =
(F \cdot H) G - H (F \cdot G)$, we have
\[
\overrightarrow{Z}= \frac{1}{ 4 \pi} \bigg{[} \partial_{y_k} \eta
\nabla_y \int_{{\bf R}^3} \frac{ G(w, l; z+y, s)}{|z|} dz -
\big{(} \sum^3_{j=1}
\partial_{y_j} \eta \partial_{y_j} \int_{{\bf
R}^3} \frac{ G(w, l; z+y, s)}{|z|} dz  \big{)} e_k \bigg{]}
\]Since $\partial_{w_i} G(w, l; y+z, s) =
-\partial_{y_i} G(w, l; y+z, s)$, the above shows
\[
\aligned
&\partial_{w_i} \overrightarrow{Z}\\
&=\frac{1}{ 4 \pi} \bigg{[}
\partial_{y_k} \eta \nabla_y \int_{{\bf R}^3} \frac{
\partial_{w_i} G(w, l; z+y, s)}{|z|} dz - \big{(} \sum^3_{j=1}
\partial_{y_j} \eta \partial_{y_j} \int_{{\bf
R}^3} \frac{\partial_{w_i} G(w, l; z+y, s)}{|z|} dz  \big{)} e_k
\bigg{]}. \endaligned \]Hence
\[
\aligned
 &\partial_{w_i} \overrightarrow{Z}
  =-\frac{1}{ 4 \pi} \bigg{[}
\partial_{y_k} \eta \nabla_y \int_{{\bf R}^3} \frac{
\partial_{y_i} G(w, l; z+y, s)}{|z|} dz \\
&\qquad \qquad - \big{(} \sum^3_{j=1}
\partial_{y_j} \eta \partial_{y_j} \int_{{\bf
R}^3} \frac{\partial_{y_i} G(w, l; z+y, s)}{|z|} dz  \big{)} e_k
\bigg{]}.
\endaligned
\leqno(3.16')
\]Substituting the above into the defining formula for $M_2$
((3.15)) and use integration by parts, we obtain
\[
\aligned M_2&= \int_{Q_r(x, t)} \partial_{y_i} u(y, s) [\Delta
\overrightarrow{Z} +
\partial_s \overrightarrow{Z}](y, s) dyds\\
&\qquad + \frac{1}{ 4 \pi} \int_{Q_r(x, t)} u(y, s) (\Delta  +
\partial_s) [ \partial_{y_i}
\partial_{y_k} \eta \nabla_y \int_{{\bf
R}^3} \frac{ G(w, l; z+y, s)}{|z|}] dz dyds\\
&\qquad + \frac{1}{ 4 \pi} \int_{Q_r(x, t)} u(y, s) (\Delta  +
\partial_s) [ \sum^3_{j=1} \partial_{y_i}
(\partial_{y_j} \eta)   \partial_{y_j} \int_{{\bf R}^3} \frac{
G(w, l; z+y, s)}{|z|} dz] dyds
\endaligned
\]Just like the estimate of $J_2$ in part (a), we have
\[
|\Delta \overrightarrow{Z} +\partial_s \overrightarrow{Z}| \le
c/r^5.
\]Here, as before, we used the fact that the parabolic
distance between $(w, l)$ and $(y, s)$ is a least $r/4$.

For the rest of the terms, using (3.12) and the same argument as
in the estimate of $K_1$, we deduce
\[
|M_2| \le \frac{c}{r^5} \int_{Q_r(x, t)-Q_{r/2}(x, t)} |\nabla
u(y, s)| dyds +\frac{c}{r^6} \int_{Q_r(x, t)-Q_{r/2}(x, t)} |u(y,
s)| dyds. \leqno(3.17)
\]

\medskip

Finally we bound $M_3$. Using integration by parts on the $y$
variable and using the assumption that $div \ b=0$, we can write
\[
\aligned M_3&=  \int_{Q_r(x, t)} u b \nabla_y
\partial_{w_i} \Phi^{(w, l)}_k(y, s) dyds \\
&= - \int_{Q_r(x, t)} \partial_{w_i} \Phi^{(w, l)}_k(y, s)  b
\nabla_y u dyds \\
&= - \int_{Q_r(x, t)} (\partial_{w_i} (\eta E_k +
\overrightarrow{Z}) )  b \nabla u dyds.
\endaligned
\]The last step is by (3.8), with $(w, l)$ replacing $(x, t)$ there. From the well known property of the Stokes system
\[
|\partial_{w_i} E_k(w, l; y, s) | \le C K_1(w, l; y, s).
\]This together with (3.16') and (3.12) yield
\[
|\partial_{w_i} (\eta E_k + \overrightarrow{Z})| \le C K_1(w, l;
y, s)
\]

Therefore
\[
|M_3| \le C \int_{Q_r(x, t)} K_1(w, l; y, s)  |b \nabla u| dyds.
\leqno(3.18)
\]Recall (see Remark 3.1 just before Lemma 3.2),  the integral on
the righthand side of (3.18) is a $L^1_{loc}$ function of space
time.

Combining (3.15)-(3.18), we have, for $(w, l) \in Q_{r/4}(x, t)$,
a.e.
\[
\aligned
| \int_{{\bf R}^3}& u(w, l) \partial_{w_i} \rho(w) dw |\\
 &\le \frac{c}{r^5} \int_{Q_r(x, t)-Q_{r/2}(x, t)} |\nabla u(y,
s)| dyds +\frac{c}{r^6}
\int_{Q_r(x, t)-Q_{r/2}(x, t)} |u(y, s)| dyds\\
&\qquad + c \int_{{\bf R}^3} \int_{Q_r(x, t)} K_1(w, l; y, s)  |b|
\ |\nabla u| dyds \ \rho(w) dw.
\endaligned
\]Since $\rho$ is arbitrary and $\nabla u$ and $M_3$ is
$L^1_{loc}$, we deduce, for $(w, l) \in Q_{r/4}(x, t)$,
\[
\aligned
|\partial_{w_i} &u(w, l)|\\
 &\le \frac{c}{r^5} \int_{Q_r(x, t)-Q_{r/2}(x, t)} |\nabla u(y,
s)| dyds +\frac{c}{r^6}
\int_{Q_r(x, t)-Q_{r/2}(x, t)} |u(y, s)| dyds\\
&\qquad +  \int_{Q_r(x, t)} K_1(w, l; y, s)  |b| \ |\nabla u|
dyds.
\endaligned
\]
This proves the lemma. \qed

\medskip

Now we are ready to give the

{\bf Proof of Theorem 1.2.}
\medskip

{\it (a) We will  prove the first formula, i.e. the bound on $u(x,
t)$.}

 The idea of the
proof is to iterate (3.6) in a special manner since simple
iteration will double the domain of integration in each step and
will not yield a local formula. The key is to cut in half the size
of the cube in (3.6) after each iteration. Here are the details.

In order to simplify the presentation, we will use capital letters
to denote points in space-time. For example $X=(x, t), Y=(y, s),
Z=(z, \tau)$, etc. From (3.6), we have
\[
|u(X)| \le \lambda \frac{1}{r^5} \int_{Q_r(X)} |u(Y)| dY + \lambda
\int_{Q_r(X)} K_1(X; Y) |b(Y)| \ |u(Y)| dY. \leqno(3.19)
\]For $Y \in Q_r(X)$, we apply the representation formula for
cubes of half of the previous size to get
\[
|u(Y)| \le \lambda \frac{2^5}{r^5} \int_{Q_{r/2}(Y)} |u(Z)| dZ +
\lambda \int_{Q_{r/2}(Y)} K_1(Y; Z) |b(Z)| \ |u(Z)| dZ.
\]Combining the above two inequalities, we obtain
\[
\aligned |u(X)| &\le \lambda \frac{1}{r^5} \int_{Q_r(X)} |u(Y)| dY
+ \lambda \int_{Q_r(X)} K_1(X; Y) |b(Y)| \ \lambda \frac{2^5}{r^5}
\int_{Q_{r/2}(Y)} |u(Z)| dZ dY \\
&\qquad + \lambda \int_{Q_r(X)} K_1(X; Y) |b(Y)| \lambda
\int_{Q_{r/2}(Y)} K_1(Y; Z) |b(Z)| |u(Z)| \  dZ dY.
\endaligned
\]Notice that $Q_{r/2}(Y) \subset Q_{3r/2}(X)$ when $Y \in
Q_r(X)$. The above thus shows
\[
\aligned |u(X)| &\le \lambda \frac{1}{r^5} \int_{Q_r(X)} |u(Y)| dY
+ \lambda^2 \frac{2^5}{r^5} \Vert u \Vert_{L^1(Q_{2r}(X))}
\int_{Q_r(X)} K_1(X; Y) |b(Y)| dY \\
&\qquad + \lambda^2 \int_{Q_{3r/2}(X)} \ \int_{Q_r(X)} K_1(X; Y)
|b(Y)|  K_1(Y; Z) dY |b(Z)| |u(Z)| dZ.
\endaligned
\]

Applying Lemma 3.1, we deduce
\[
\aligned
 |u(X)| &\le \lambda \frac{1}{r^5}  \Vert u
\Vert_{L^1(Q_{2r}(X))} + \lambda^2 \frac{2^5}{r^5} \Vert u
\Vert_{L^1(Q_{2r}(X))} B(b) \\
&\qquad + \lambda^2  c B(b) \int_{Q_{3r/2}(X)} K_1(X; Z) |b(Z)|
|u(Z)| dZ,
\endaligned
\]here and throughout this section
$ B(b)=B(b, t-(4r)^2, t)$ defined in (1.6').

For $u(Z)$, we use the representation formula in the cube
$Q_{r/4}(Z)$:
\[
|u(Z)| \le \lambda \frac{4^5}{r^5} \int_{Q_{r/4}(Z)} |u(W)| dW +
\lambda \int_{Q_{r/4}(W)} K_1(Z; W) |b(W)| \ |u(W)| dW.
\]Note that for $Z \in Q_{3r/2}(X)$,  we have
$Q_{r/4}(Z) \subset Q_{7r/4}(X) \subset Q_{2r}(X)$. Hence
\[
\aligned |u(X)| &\le \lambda \frac{1}{r^5}  \Vert u
\Vert_{L^1(Q_{2r}(X))} +  \lambda^2 \frac{2^5}{r^5} \Vert u
\Vert_{L^1(Q_{2r}(X))} B(b)  \\
&+ \lambda^2 \frac{4^5}{r^5} c B(b)
\int_{Q_{3r/2}(X)}  K_1(X; Z) |b(Z)|dZ \,  \Vert u \Vert_{L^1(Q_{2r}(X))} \\
&+ \lambda^3  c B(b) \int_{Q_{3r/2}(X)}  K_1(X; Z) |b(Z)|
\int_{Q_{r/4}(Z)} K_1(Z; W) |b(W)| |u(W)| dW dZ.
\endaligned
\]Exchanging the integrals and using Lemma 3.1 again, we deduce
\[
\aligned |u(X)| &\le \lambda \frac{1}{r^5}  \Vert u
\Vert_{L^1(Q_{2r}(X))} + \lambda^2 \frac{2^5}{r^5} \Vert u
\Vert_{L^1(Q_{2r}(X))} B(b) + \lambda^2 \frac{4^5}{r^5} c B(b)^2 \Vert u \Vert_{L^1(Q_{2r}(X))}
 \\
&+ \lambda^3  c^2 B(b)^2 \int_{Q_{7r/4}(X)}  K_1(X; W)  |b(W)|
|u(W)| dW dZ.
\endaligned
\]

We iterate the above process, halving the size of each cube. By
induction, it is clear that for some $C, c_1>0$,
\[
|u(X)| \le \frac{C}{r^5}  \Vert u \Vert_{L^1(Q_{2r}(X))}
\Sigma^\infty_{j=1} [2^5 c_1  B(b)]^j.
\]When $2^5 c_1 B(b)<1$ the above series converges to yield the mean
value inequality for $u$. Since $b$ is in the Kato class defined
in Definition 1.2, we know that $2^5 c_1 B(b)= 2^5 c_1 B(b,
t-(4r)^2, t)<1$ when $r$ is sufficiently small. This proves the
bound on $u$.

\medskip

{\it (b) We prove the gradient bound.}
\medskip

Since $u+c$ is also a solution of (1.2) for any constant $c$, we
will assume that $\overline{u}_{Q_{2r}}$, the average of $u$ in
$Q_{2r}(x, t)$ is $0$.

The idea is to iterate (3.6') in the above manner. From (3.6'),
using the same notation as in part (a), we have
\[
 |\nabla u(X)| \le  m(X, r)
 + \lambda \int_{Q_r(X)} K_1(X; Y) |b(Y)| \ |\nabla
u(Y)| dY, \leqno(3.20)
\]where
\[
m(X, r) \equiv
 \frac{\lambda}{r^5} \int_{Q_r(X)} | \nabla u(Y)|  dyds
 + \frac{\lambda}{r^6} \int_{Q_r(X)} |u(Y)| dY.
\]Here, we remind the reader that both sides of (3.20) are
$L^1_{loc}$ functions and hence finite almost everywhere (see
Remark 3.1 just before Lemma 3.2).

Applying (3.20) to $\nabla u(Y)$ and $Q_{r/2}(Y)$, we obtain
\[
|\nabla u(Y)| \le  m(Y, r/2)
 + \lambda \int_{Q_{r/2}(Y)} K_1(Y; Z) |b(Z)| \ |\nabla
u(Z)| dZ.
\]For $Y \in Q_r(X)$, it is clear that there exists a $\mu >0$,
independent of $r$ such that
\[
m(Y, r/2) \le \mu m(X, 2 r).
\]Hence
\[
|\nabla u(Y)| \le  \mu m(X, 2 r)
 + \lambda \int_{Q_{r/2}(Y)} K_1(Y; Z) |b(Z)| \ |\nabla
u(Z)| dZ. \leqno(3.21)
\]Substituting (3.21) to (3.20) we have
\[
\aligned |\nabla u(X)| &\le  m(X, r) + \mu \ m(X, 2r) \lambda
\int_{Q_r(X)} K_1(X; Y) |b(Y)| dy \\
&\qquad + \lambda \int_{Q_r(X)} K_1(X; Y) |b(Y)| \lambda
\int_{Q_{r/2}(Y)} K_1(Y; Z) |b(Z)| \ |\nabla u(Z)| dZ dY.
\endaligned
\]Therefore
\[
\aligned |\nabla u(X)| &\le  m(X, r) + \mu \lambda \ m(X, 2r) B(b) \\
&\qquad + \lambda^2 \int_{Q_{3r/2}(X)} \int_{Q_{r/2}(Z)} K_1(X; Y)
|b(Y)| K_1(Y; Z) dZ |b(Z)| \ |\nabla u(Z)| dZ dY.
\endaligned
\]Lemma 3.1 then implies
\[
|\nabla u(X)| \le  m(X, r) + \mu \lambda \ m(X, 2r) B(b) \\
 + c \lambda^2 B(b) \int_{Q_{3r/2}(X)} K_1(X,Z) |b(Z)| \ |\nabla
u(Z)| dZ.
\]Now, using (3.20) on $|\nabla
u(Z)|$ and the cube $Q_{r/4}(Z)$ and repeat the above argument,
halving the size of the cube in each step, we finally reach
\[
|\nabla u(X)| \le  C  \ m(X, 2r) \Sigma^{\infty}_{k=1}
\lambda^{k+1} \mu^{k+1} B(b)^k.
\]As before this implies the desired gradient bound when $B(b)$ is
small.

The last statement of Theorem 1.2 is a simple consequence of the
gradient estimate and Poincar\`e inequality. \qed

We end the section by showing that if the sum of the entries of
the drift term $b$ is zero, then (1.2) in torus has bounded
solutions for many initial values, regardless of the singularity
of $b$. To state the result rigourously, we will assume that $b$
is bounded. However all coefficients are independent of the bounds
of $b$.

\medskip

\begin{proposition} Given bounded vector fields $b=(b_1(x, t), ...,
b_n(x, t))$, consider the linearized Navier-Stokes equation with
periodic boundary condition i.e. in a torus.
\[
\begin{cases}
\Delta u(x, t) - b(x, t)  \nabla u(x, t) + \nabla P(x, t) -
\partial_t
u(x, t)=0, \quad (x, t) \in D  \times {\bf R},\\
div u = 0, \ div b=0, \ b(\cdot, t) \in L^\infty \\
u(x, 0) = u_0(x).
\end{cases}
\leqno(3.22)
\]Here $D = [0, 2 \pi]^n$. The functions $b(\cdot, t)$,
$u_0(\cdot)$ and $u(\cdot, t)$ have period $2 \pi$.

Suppose $\sum^n_{j=1} b_j(x, t) = \lambda$, a constant and $u_0$
is any finite linear combination of
\[
\int_{D} e^{i k \sum^n_{j=1}(x_j-y_j)} f(y) dy,
\]where $k$ is a positive integer and $f$ is a bounded,
divergence free vector field with period $2 \pi$. Then there
exists a constant $c$ independent of $b$ and $\lambda$ such that
\[
u(x, t) \le c e^{-t} C(\Vert u_0 \Vert_{\infty}).
\]
\proof
\end{proposition}

Under the assumption that $b$ is bounded, the existence of
solutions to (3.22) follows from the standard theory. Let $E=E(x,
t; y, s)$ be the fundamental solution of (3.22). The existence of
$E$ is also standard.

 First we
just assume that $u_0$ has  one term, i.e.
\[
u_0(x) = \int_{D} e^{i k \sum^n_{j=1}(x_j-y_j)} f(y) dy,
\]Fixing $(x, t)$, consider
\[
I(s) \equiv \int_D E(x, t; y, s) u_0(y) dy. \leqno(3.23)
\]By the fact that the rows of $E$ satisfies the conjugate
equation of (3.22), we have
\[
\aligned I'(s) & = \int_{D} \frac{d}{ds} E(x, t; y, s) u_0(y) dy\\
&= \int_{D} \big{[} - \Delta_y E(x, t; y, s) - b(y, s) \nabla_y
E(x, t; y, s) + \nabla Y(y, s)] u_0(y) dy. \endaligned
\]Here
\[
\nabla Y =\left( \begin{array}{cc} \partial_1 P_1 \ ...\
\partial_n P_1\\... \ ... \ ... \\ \partial_1 P_1\ ... \
\partial_n P_n
\end{array} \right)
\]with $P_1, ..., P_n$ being scalar functions. Using integration
by parts and the divergence free property of $b$, we deduce
\[
I'(s) = - \int_{D}   E(x, t; y, s) \Delta u_0(y) dy  + \int_{D}
E(x, t; y, s) \sum^n_{j=1} b_j \partial_{y_j} u_0(y) dy.
\]Noticing that $\Delta u_0 = - n k^2 u_0$ and $\partial_{y_j} u_0(y)
= i k u_0$, we have
\[
I'(s) =  n k^2 \int_{D}   E(x, t; y, s) u_0(y) dy + i k \int_{D}
(\sum^n_{j=1} b_j) E(x, t; y, s)  u_0(y) dy.
\]By our assumption that $\sum^n_{j=1} b_j = \lambda$, the above
shows
\[
I'(s) =  (n k^2 + i k \lambda) I(s).
\]Hence
\[
I(s) = e^{(n k^2 + i k \lambda) s} I(0).
\]From (3.23) and the fact that $E$ is the fundamental solution to
(3.22), we have $I(0) = u(x, t)$ and $I(t) = u_0(x)$. This shows
\[
u(x, t) = e^{- (n k^2 + i k \lambda) t} u_0(x). \leqno(3.24)
\]

Now let $S$ be a set of finite positive integers and
\[
u_0(x) = \sum_{l \in S} c_l \int_{D} e^{i k_l
\sum^n_{j=1}(x_j-y_j)} f_l(y) dy.
\]Here $k_l$ is a positive integer and $f_l$ is a bounded,
divergence free vector field with period $2 \pi$. By (3.24) one
has
\[
u(x, t) = \sum_{l \in S} e^{- (n k^2_l + i k_l \lambda) t} c_l
\int_{D} e^{i k_l \sum^n_{j=1}(x_j-y_j)} f_l(y) dy.
\]\qed

\section{A regularity condition for Navier-Stokes equations not involving
absolute values}
%sec5

In this section we introduce another sufficient condition on the
velocity for boundedness of weak solutions of 3 dimensional
Navier-Stokes equations. The novelty is that no absolute value of
$u$ is involved. This is useful since it allows more cancellation
effect to be taken into account. Throughout the years, various
conditions on $u$ that imply regularity have been proposed. One of
them is the Prodi-Serrin condition which requires that $u \in
L^{p, q}$ with $\frac{3}{p}+\frac{2}{q} \le 1$ for some $3<p \le
\infty$ and $q \ge 2$. See ([Ser], [Str] e.g.) Recently the
authors in [ESS] showed that the condition $p=3$ and $q=\infty$
also implies regularity. In another development the author of [Mo]
improved the Prodi-Serrin condition by a log factor, i.e. by
requiring
\[
\int^T_0 \frac{\Vert u(., t) \Vert^q_p}{1+ \log^+ \Vert u(.,t)
\Vert_p} dt < \infty,
\]where $3/p + 2/q = 1$ and $3<p<\infty$, $2<q<\infty$

Most recently in [Z3], a form boundedness condition on velocity
was introduced, which will imply the boundedness of weak
solutions. The form boundedness condition,  with its root in the
perturbation theory of elliptic operators and mathematical
physics, seems to be different from all the previous conditions.
It seems to be one of the most general condition under the
available tools. This fact has been well documented in the theory
of linear elliptic equations. See [Si] e.g.  Moreover, as
indicated in [Z3],  it contains the Prodi-Serrin condition except
when $p$ or $q$ are infinite. It also includes suitable
Morrey-Compamato type spaces. However we are not sure this
condition contains the one in [Mo].

More precisely we proved
\medskip

\noindent {\it {\bf Theorem 5.1 ([Z3])}
%Thm4.1
Let $u$ be a Leray-Hopf solution to the 3 dimensional
Navier-Stokes equation  in ${\bf R}^3 \times (0, \infty)$.
\[
\begin{cases}
\Delta u(x, t) - u(x, t)  \nabla u(x, t) + \nabla P(x, t) -
\partial_t
u(x, t)=0, \quad (x, t) \in \Omega \subset {\bf R}^3 \times {\bf R},\\
div u = 0.
\end{cases}
\leqno(4.1)
\]
 Suppose
for every $(x_0, t_0) \in {\bf R}^3 \times (0, \infty)$, there
exists a cube $Q_r=B(x_0, r) \times [t_0-r^2, t_0]$ such that  $u$
satisfies the form bounded condition
\[
\int_{Q_r} |u|^2 \phi^2 dyds \le \frac{1}{24} \bigg{(} \int_{Q_r}
|\nabla \phi|^2 dyds + \sup_{s \in [t_0-r^2, t_0]} \int_{B(x_0,
r)} \phi^2(y, s) dy \bigg{)} + B(\Vert \phi \Vert_{L^2(Q_r)}).
\leqno(4.2)
\]Here $\phi$ is any smooth function vanishing on the parabolic side of
$Q_r$ and $B=B(t)$ is any given locally bounded function of $t \in
{\bf R}^1$. Then $u$ is a classical solution when $t>0$.}
\medskip

In this section we are able to extend the form boundedness
condition further. The main improvement is that our new condition
in $u$ ((4.3)) below does not involve the absolute value of $u$.
This differs significantly from known conditions on $u$ so far
where $|u|$ is always present. By a simple integration by parts
argument, it is clear that Condition (4.3) is more general than
Condition (4.2).

\medskip

\begin{theorem}
%Thm4.1
Let $u$ be a Leray-Hopf solution to the 3 dimensional
Navier-Stokes equation  in ${\bf R}^3 \times (0, \infty)$.

 Suppose
for every $(x_0, t_0) \in {\bf R}^3 \times (0, \infty)$, there
exists a cube $Q_r=B(x_0, r) \times [t_0-r^2, t_0]$ such that  $u$
satisfies the form bounded condition: for a given $\delta>0$,
\[
\int_{Q_r} \phi \nabla u \cdot \phi dyds \le \frac{1-\delta}{2}
\bigg{(} \int_{Q_r} |\nabla \phi|^2 dyds + \frac{1}{2} \sup_{s \in
[t_0-r^2, t_0]} \int_{B(x_0, r)} \phi^2(y, s) dy \bigg{)} +
B(\Vert \phi \Vert_{L^2(Q_r)}). \leqno(4.3)
\]Here $\phi$ is any smooth vector field vanishing on parabolic the side of
$Q_r$ and $B=B(t)$ is any given locally bounded function of $t \in
{\bf R}^1$. Then $u$ is a classical solution when $t>0$.
\end{theorem}
\medskip

{\bf Remark 4.1.} Condition 4.3 is actually a condition on the
strain tensor $\nabla u + (\nabla u)^T$. Theorem 4.1 immediately
implies that weak solutions to the 3 dimensional Navier-Stokes
equations are locally bounded in any open subset of the region
where the strain tensor is negative definite.
\medskip

{\bf Proof of Theorem 4.1.} Let $t_0$ be the first moment of
singularity formation. We will reach a contradiction. It is clear
that we only need to prove that $u$ is bounded in
$Q_{r/8}=Q_{r/8}(x_0, t_0)$ for some $r>0$. In fact the number $8$
is not essential. Any number greater than $1$ would work.

 Consider the equation for
vorticity $w = \nabla \times u$. It is well known that, in the
interior of $Q_r$,  $w$ is a classical solution to the parabolic
system with singular coefficients
\[
\Delta w - u \nabla w + w \nabla u - w_t = 0. \leqno(4.4)
\]Let $\psi=\psi(y, s)$ be the refined cut-off function defined right after
(2.1) such that $\psi=1$ in $Q_{r/2}$, $\psi=0$ in $Q^c_r$ and
such that $0 \le \psi \le 1$, $|\nabla \psi| \le C/r$ and
$|\psi_t| \le C/r^2$. We can use $w \psi^2$ as a test function on
(4.4) to obtain
\[
\aligned
 \int_{Q_r} | &\nabla (w \psi) |^2 dyds + \frac{1}{2}
 \int_{B(x_0, r)} |w \psi|^2(y, t_0) dy\\
&\le \frac{C}{r^2} \int_{Q_r} |w|^2 dyds - \int_{Q_r} u \nabla w
\cdot w \psi^2 dyds  + \int_{Q_r} w \nabla u \cdot w \psi^2
dyds \\
&\equiv I_1 + I_2 + I_3.
\endaligned
\leqno(4.5)
\]

The term $I_1$ is already in good shape. Next, using integration
by parts and the divergence free condition on $u$, we have
\[
I_2 = \frac{1}{2} \int_{Q_r} u \cdot \nabla \psi \psi |w|^2 dyds.
\]Since $\nabla u \in L^2_{loc}({\bf R}^3 \times [0, \infty))$
and $\Vert u( \cdot, t)\Vert_{{\bf R}^3}$ is non-increasing in
time, it is easy to prove by Sobolev imbedding and H\"older's
inequality that $u|_{Q_r}$ satisfies the form boundedness
condition (1.3). In fact this has been proven in Corollary 2 in
[Z2]. Hence we can bound $I_2$ in exactly the same way as $T_2$ in
(2.4) with $b$ being chosen as $u$ here. Following the argument
between (2.4) and (2.6), we obtain, for any given $\delta>0$,
\[
I_2 \le  \frac{\delta}{4} \int_{Q_{r}} |\nabla (\psi  w)|^2 dyds +
c_\delta e^{( c_1/ r)^2} \int_{Q_r}  w^2 dyds. \leqno(4.6)
\]Note that in (2.6), $\delta$ was chosen as $2$. However, a
closer look at the proof shows that (4.6) is true.

Next we estimate $I_3$. From Condition (4.3)
\[
I_3 \le \frac{1-\delta}{2} \bigg{(}  \int_{Q_r} |\nabla (w \psi)
|^2 dyds + \frac{1}{2} \sup_{s \in [t_0-r^2, t_0]} \int_{B(x_0,
r)} (w \psi)^2(y, s) dy \bigg{)} + B(\Vert w \psi
\Vert_{L^2(Q_r)}). \leqno(4.7)
\]

Substituting (4.6)-(4.7) to (4.5) we obtain,
\[
\aligned
 \int_{Q_r}& | \nabla (w \psi) |^2 dyds + \frac{1}{2}
\int_{B(x_0, r)} |w \psi|^2(y, t_0) dy \\
&\le \frac{\delta}{4} \int_{Q_{r}} |\nabla (\psi  w)|^2 dyds +
c_\delta e^{( c_1/ r)^2}
\int_{Q_r}  w^2 dyds \\
&+ \frac{1-\delta}{2} \bigg{(}  \int_{Q_r} |\nabla (w \psi) |^2
dyds + \frac{1}{2} \sup_{s \in [t_0-r^2, t_0]} \int_{B(x_0, r)} (w
\psi)^2(y, s) dy \bigg{)}
 + C B(\Vert w
\Vert_{L^2(Q_r)}).
\endaligned
\]Repeating the above process, but restrict the integrals to $Q_r
\cap \{ (y, s) \ | s<t \}$ with $t<t_0$, we obtain, for any $t \in
[t_0-r^2, t_0]$,
\[
\aligned &\frac{1}{2}
\int_{B(x_0, r)} |w \psi|^2(y, t) dy \\
&\le \frac{\delta}{4} \int_{Q_{r}} |\nabla (\psi  w)|^2 dyds +
c_\delta e^{( c_1/ r)^2}
\int_{Q_r}  w^2 dyds \\
&+ \frac{1-\delta}{2} \bigg{(}  \int_{Q_r} |\nabla (w \psi) |^2
dyds + \frac{1}{2} \sup_{s \in [t_0-r^2, t_0]} \int_{B(x_0, r)} (w
\psi)^2(y, s) dy \bigg{)}
 + C B(\Vert w
\Vert_{L^2(Q_r)}).
\endaligned
\]

Combining the last two estimates, we deduce,
\[
\aligned
 &\int_{Q_r} | \nabla (w \psi) |^2 dyds + sup_{t_0-r^2 \le
s \le t_0} \int_{B(x_0, r)} |w \psi|^2(y, s) dy \\
&\le \frac{C_1}{r^2} \Vert w \Vert_{L^2(Q_r)} + B_1(\Vert w
\Vert_{L^2(Q_r)}). \endaligned
 \leqno(4.8)
\]Here $B_1$ is a locally bounded, one variable function.

Using standard results, we know that (4.8) implies that $u$ is
regular. Here is the proof.

From (4.8), it is clear that $\int_{Q_r} | curl \ (w \psi) |^2
dyds \le C.$ Hence, since $\psi =1$ in $Q_{r/2}$,
\[
\int_{Q_{r/2}} | \Delta  u|^2 dyds \le C. \leqno(4.9)
\]Let $\eta =\eta(y)$ be a cut-off function such that $\eta=1$ in
$B(x_0, r/4)$ and $\eta=0$ in $B(x_0, r/2)^c$. Then for each $s
\in [t_0-(r/4)^2, t_0]$, we have, in the weak sense
\[
\Delta (u \eta) = \eta \Delta u + 2 \nabla u \nabla \eta + u
\Delta \eta \equiv f,
\]in $Q_{r/2}$. By standard elliptic estimates, using the fact
that $u \eta=0$ on the boundary,
\[
\Vert D^2 u( \cdot, s) \Vert_{L^2(B(x_0, r/4))} \le C \Vert f(
\cdot, s) \Vert_{L^2(B(x_0, r/2))}.
\]This shows that
\[
\Vert D^2 u \Vert_{L^2(Q_{r/4})} \le C \Vert \Delta u
\Vert_{L^2(Q_{r/2})} + C \Vert u \Vert_{L^2(Q_{r/2})}.
\]By Sobolev imbedding
\[
\nabla u \in L^{6, 2}(Q_{r/4}). \leqno(4.10)
\]

Next, from (1.22) on p316 [Te],
\[
\Vert u( \cdot, s) \eta \Vert_{W^{1, 2}} \le C \big{(} \Vert u(
\cdot, s) \eta \Vert_{L^2} + \Vert div (u \eta) ( \cdot, s)
\Vert_{L^2} + \Vert curl (u( \cdot, s) \eta) \Vert_{L^2} \big{)}.
\]Here all norms are over the ball $B(x_0, r/2)$. Therefore
\[
\Vert u \eta ( \cdot, s)  \Vert_{W^{1, 2}} \le C \big{(} \Vert u(
\cdot, s) \eta \Vert_{L^2} + \Vert u \nabla \eta ( \cdot, s)
\Vert_{L^2} + \Vert w \eta ( \cdot, s)  \Vert_{L^2} + \Vert |u(
\cdot, s)| \ |\nabla \eta| \Vert_{L^2} \big{)}.
\]It follows
\[
\Vert u( \cdot, s)  \Vert_{W^{1, 2}(B(x_0, r/4))} \le C.
\]From Sobolev imbedding we know that
\[
u \in L^{6, \infty}(Q_{r/4}). \leqno(4.11)
\]We treat $u$ and $\nabla u$ as coefficients in the vortex equation (4.4).
By (4.10) and (4.11), the standard parabolic theory (see [Lieb]
e.g.), shows that $w$ is bounded and  H\"older continuous in
$Q_{r/8}$. Here the bound depends only on the $L^2$ norm of $w$ in
$Q_r$ and $r$. This is so because of the relation $ 3/6 + 2/\infty
< 1$ for the norm of $u$ and $3/6 + 2/2 < 2$ for the norm of
$\nabla u$. Now a standard bootstrapping argument shows that $u$
is smooth. \qed

\medskip

\medskip

{\bf Acknowledgement} I thank Professors Z. Grujic, I. Kukavica
and V. Liskevich for helpful discussions.

\medskip

\section{Addendum May 2019} 

The purpose of the addendum is to fill in a missing term in Theorem 1.7 and Lemma 3.3 in the current paper that is published in Pacific Journal of Mathenatics, Vol. 223, No.2, 2006, which will be cited as \cite{Z:4}.  This is due to the use of a formula in a cited reference, which omitted a term.
The main conclusion that local solutions of certain linearized Navier-Stokes equation have bounded spatial gradient is intact. This includes bounded local Leray-Hopf solutions to the Navier Stokes equation without condition on the pressure.

Let us recall the basic set up of \cite{Z:4}. The equation is: for $(x, t) \in \Omega \subset {\bf R}^3 \times {\bf R}$,
\be
\lab{lns}
\begin{cases}
\Delta u(x, t) - b(x, t)  \nabla u(x, t) + \nabla P(x, t) -
\partial_t
u(x, t)=0, \\
div \, u = 0, \ div \, b=0, \ b(\cdot, t) \in L^2_{loc}.
\end{cases}
\ee
Here $\Delta$ is the standard Laplacian and $b=b(x, t)$ is a  given
$L^2_{loc}$  vector field  to be specified later. $\Omega$
is a domain.
The parabolic Kato norm is:
\be
B(b, l, t) = \sup_{x} \int^{t}_l \int_{{\bf R}^n} [K_1(x, t; y, s) +
K_1(x, s; y, l) ] |b(y, s)| dyds.
\lab{(1.6')}
\ee We say $b$ is in class $K_1$ if $\lim_{l \to t} B(b, l, t)=0$. $K_1$ class in a bounded space time domain contains the usual $L^{p, q}$ functions with $(3/p)+(2/q)<1.$

Let $D$ for a domain in ${\bf R}^3$ and $T>0$, we will use this definition for solutions of (\ref{lns}).

\begin{definition}
\lab{defsol}
 A divergence free vector
field $u \in L^{\infty}(0, T; L^2(D)) \cap L^2(0, T; W^{1, 2}(D))$
is called a
 (weak) solution of (\ref{lns}) if:
for any vector valued $\phi \in C^{\infty}(D \times [0, T])$ with
$ div \, \phi =0$ and $\phi=0$ on $\partial D \times [0, T]$, $u$
satisfies
\[
\int^{t_2}_{t_1}\int_{{\bf R}^n} <u, \partial_t \phi + \Delta
\phi>
 dxdt - \int^{t_2}_{t_1}\int_{{\bf R}^n} <b \nabla u,  \phi>  dxdt
 = - \int_{{\bf R}^n} <u, \phi>(x, t) |^{t_2}_{t_1} dx.
\]
\end{definition}

The corrected version of Theorem 1.7 in \cite{Z:4} is:

\begin{theorem}
\label{th:1.2}
Let $u$ be a solution of (\ref{lns}) in a domain $\Omega \subset {\bf
R}^3 \times {\bf R}$.  Suppose $Q_{4r}(x, t) \subset \Omega$, $div
b =0$ and that $b |_{Q_{2r}(x, t)}$ is in class $K_1$ and $b \in
L^2_{loc}$. Then both $u$ and $|\nabla u|$ are bounded functions
in $Q_{2r}=Q_{2r}(x, t)$, the standard parabolic cube of size $2r$.

Moreover, for some positive constants $C=C(b)$ and $r_0$,
depending on the size of the Kato norm of $b$,   there hold, when
$0<r<r_0$,
\[
\aligned
 |u(x, t)| &\le  \frac{C}{r^3} \sup_{s \in [t-(2 r)^2, t]} \int_{B_{r}(x)} |u(y,
s)| dy,\\
|\nabla u(x, t)| &\le \frac{C}{r^5} \Vert \nabla
u \Vert_{L^1(Q_{2r})}
 + \frac{C}{r^4} \sup_{s \in [t-(2 r)^2, t]} \int_{B_{r}(x)} |u(y,
s)| dy.
\endaligned
\]
\end{theorem}

The proof of the theorem is through a scaled iteration process starting with Lemma 3.3 in \cite{Z:4}, which contains mean value inequalities for: (a) solutions of
(\ref{lns}); (b)  spatial gradient of solutions. Inequality (a) is first obtained in \cite{O:1}, which is also the basis for
inequality (b).  Recently Professor Hongjie Dong kindly informed us that inequality (a) and hence (b) misses a term.  The corrected version of Lemma 3.3 in \cite{Z:4} is the following. In our opinion the main localization idea in \cite{O:1} is still nice.

\begin{lemma} (mean value inequalities, replacing Lemma 3.3 in \cite{Z:4})
%L3.2.

(a). Let $u$ be a solution of (\ref{lns}) in the region $\Omega$.
Suppose $Q_{2r}(x, t) \subset \Omega$. Then there exists a
constant $\lambda$ such that
\[
\aligned
 |u(x, t)| &\le  \frac{\lambda}{r^5} \int_{Q_r(x,
t)-Q_{r/2}(x, t)} |u(y, s)| dyds + \underbrace{\frac{\lambda}{r^3} \int_{B_r(x)-B_{r/2}(x)} |u(y, t)| dy}\\
&\qquad + \lambda \int_{Q_r(x, t)} K_1(x, t; y, s) |b(y, s)| \
|u(y, s)| dyds. \endaligned
\]

(b). Under the same assumption as (a), there exists a constant
$\lambda$ such that
\[
%\lab{(3.6')}
\aligned
&|\nabla u(x, t)| \le  \frac{\lambda}{r^5} \int_{Q_r(x,
t)-Q_{r/2}(x, t)} | \nabla u(y, s)|  dyds
 + \frac{\lambda}{r^6}
\int_{Q_r(x, t)-Q_{r/2}(x, t)} |u(y, s)| dyds \\
&\quad + \underbrace{\frac{\lambda}{r^4} \int_{B_r(x)-B_{r/2}(x)} |u(y, t)| dy}
 + \lambda \int_{Q_r(x, t)} K_1(x, t; y, s) |b(y, s)| \
|\nabla u(y, s)| dyds. \endaligned
\]
\end{lemma}

\proof of (a). The under-braced terms were missing. We indicate the changes needed for the proof,  dividing into 3 steps.

{\it Step 1.}
Let $E=E(x, t; y, s)$ be the fundamental solution (matrix) of the
Stokes system in ${\bf R}^3 \times (0, \infty)$ and $E_k$ be the
$k$-th column of $E$. Then
\be
\lab{Ekform}
E_k(x, t; y, s) =  G(x, t; y, s) e_k + \frac{1}{4 \pi} \nabla
\partial_k \int_{{\bf R}^3} \frac{G(x, t; z, s)}{|y-z|} dz
\ee where $G$ is the fundamental solution of the heat equation.   Fixing $(x, t)$, we construct a standard
cut-off function $\eta$ such that $\eta(y, s)=1$ in $Q_{r/2}(x,
t)$, $\eta(y, s)=0$ outside of $Q_r(x, t)$, $0 \le \eta \le 1$ and
$|\nabla \eta|^2 + |\Delta \eta| + |\partial_s \eta| \le c/r^2$.

Define for $k=1, 2, 3$, after \cite{O:1},  a vector valued function
\be
\lab{(3.7)phik}
\Phi_k=\Phi^{(x, t)}_k(y, s) = \frac{1}{4 \pi} curl \big{(} \eta
(y, s) \int_{{\bf R}^3} \frac{ curl E_k(x, t; z, s)}{|z-y|} dz
\big{)}.
\ee It is clear that when $t>s$,  $\Phi_k$ is a valid test function
for equation (\ref{lns}) since $\Phi_k$ is smooth, compactly supported
and divergence free. Using $\Phi_k$ as a test function on (\ref{lns}),
by Definition \ref{defsol}, we obtain
\be
\lab{uphi1}
\int^t_0\int u(y, s) \big{(} \Delta \Phi_k + b \nabla \Phi_k +
\partial_s \Phi_k \big{)} dyds = \lim_{s \to t} \int u(y, s)
\Phi_k(y, s) dy.
\ee Here,we will suppress the superscript $(x, t)$ on
$\Phi_k$, unless there is a confusion. Also $u$ is regarded as a row vector so that $u \Phi_k$ etc is
a scalar.

Since $E_k$ is divergence free $curl \ curl E_k = - \Delta E_k$.
Thus
\be
\lab{(3.8)}
\aligned \Phi_k(y, s) &= \eta(y, s) E_k(x, t; y, s) + \frac{1}{4
\pi} \nabla \eta(y, s) \times \int_{{\bf R}^3} \frac{ curl E_k(x,
t; z, s)}{|z-y|} dz\\
&\equiv \eta E_k + \overrightarrow{Z}.
\endaligned
\ee

{\it Step 2.}
Using formula (\ref{Ekform}) of the fundamental matrix $E$ and
that of $\overrightarrow{Z}$ in (\ref{(3.8)}), direct computation shows
that
\be
\lab{uphi2}
\al
&\lim_{s \to t} \int u(y, s) \Phi_k(y, s) dy = \eta  u_k(x,t)\\
& + \int \left[
\partial_{y_k} \Gamma (x, y) (\nabla \eta \cdot u)(y, t)
+\partial_{y_k} \eta (y, t) (\nabla_y \Gamma(x, y) \cdot u(y, t) )
- (\nabla \eta \cdot \nabla_y \Gamma(x, y)) u_k(y, t) \right] dy.
\eal
\ee
where $u_k$ is the k-th component of $u$ and $\Gamma=\frac{1}{4 \pi |x-y|}$
is the Green's function on ${\bf R}^3$.

Alternatively, from (\ref{(3.8)}),
\be
\lab{i+ii}
\al
\lim_{s \to t} \int u(y, s) \Phi_k(y, s)  dy
&=\lim_{s \to t} \int \eta u(y, s) E_k(x, t; y, s) dy + \lim_{s \to t} \int u(y, s) \vec{Z} dy\\
&\equiv I + II.
\eal
\ee First we treat $I$. Consider the Helmholtz decomposition
\[
\eta u = \nabla f + X,
\]where $\Delta f = 0$, $ div \, X = 0$. Since $\Delta f = \nabla \eta \cdot u$ and
$\eta(\cdot, s)$ is compactly supported, we can take
$
f(y, s) = -\int \Gamma(y, z) ( \nabla \eta \cdot u)(z, s) dz.
$ Then
\be
\lab{x=}
X=\eta u - \nabla f = \eta u(y, s) + \int \nabla_y \Gamma(y, z) (\nabla \eta \cdot u)(z, s) dz.
\ee
Using the decay estimates
$
|\nabla_y E_k(x, t; y, s)| \le \frac{c}{(|x-y| + \sqrt{t-s})^4},
 $ $|f(y, s)| \le \frac{c}{|y|}$ and the fact that $E_k$ is a divergence vector field of variable $y$, we can integrate by parts to deduce that
$
\int \nabla f E_k(x, t; y, s) dy=0.
$ Therefore
\be
\lab{i=}
\al
I&=\lim_{s \to t} \int (\nabla f + X) E_k(x, t; y, s) dy = X_k(x, t)\\
&= \eta u_k(x, t) + \int \partial_{y_k} \Gamma(x, y) (\nabla \eta \cdot u)(y, t) dy.
\eal
\ee

In order to compute $II$, we see from (\ref{Ekform}) and (\ref{(3.8)}) that
\[
\vec{Z}
= \frac{1}{ 4 \pi} \nabla_y \eta(y, s) \times
\bigg{(} \nabla_y \big{[}\int_{{\bf R}^3} \frac{ G(x, t; z,
s)}{|z-y|} dz \big{]} \times e_k \bigg{)}.
\]This implies
\be
\lab{ii=}
\al
II=\lim_{s \to t} \int u(y, s) \vec{Z} dy
=\int  u(y, t) \cdot \left[ \nabla_y \eta \times (\nabla_y \Gamma(x, y) \times e_k )\right]  dy.
\eal
\ee

A combination of (\ref{i=}), (\ref{ii=}) and (\ref{i+ii}) proves (\ref{uphi2}). By (\ref{uphi2}) and (\ref{uphi1}):

\[
\aligned
&u_k(x, t) = \int_{Q_r(x, t)} u [ \Delta (\eta E_k) +
\partial_s (\eta
E_k)]dyds  + \int_{Q_r(x, t)} u [\Delta \overrightarrow{Z} +
\partial_s \overrightarrow{Z}] dyds\\
& + \int_{B_r(x)} \left[(\nabla \eta \cdot \nabla_y \Gamma) u_k -
\partial_{y_k} \Gamma  (\nabla \eta \cdot u)
-\partial_{y_k} \eta  (\nabla_y \Gamma \cdot u )
 \right](y, t) dy  +\int_{Q_r(x, t)} u b \nabla \Phi_k dyds.
\endaligned
\]Here and later $u b \nabla \Phi_k = u \cdot \sum^3_{i=1} b_i
\partial_i \Phi_k$.

Note that $\Delta E_k + \partial_s E_k=0$ when $s<t$. The above
infers, with $\Gamma=\Gamma(x, y)$, a local representation formula which may be of independent interest:
\be
\aligned
&u_k(x, t)
= \int_{Q_r(x, t)} u(y, s) [ E_k(x, t; y, s) (\Delta \eta +
\partial_s \eta)(y, s) +
2 \nabla \eta(y, s) \nabla_y E_k(x, t; y, s)]dyds \\
&\qquad + \int_{Q_r(x, t)} u(y, s) [\Delta \overrightarrow{Z} +
\partial_s \overrightarrow{Z}](y, s) dyds + \int_{Q_r(x, t)} u b \nabla_y \Phi_k(y, s) dyds\\
&\qquad + \int_{B_r(x)} \left[(\nabla \eta \cdot \nabla_y \Gamma) u_k -
\partial_{y_k} \Gamma  (\nabla \eta \cdot u)
-\partial_{y_k} \eta  (\nabla_y \Gamma \cdot u )
 \right](y, t) dy \\
&\equiv J_1+J_2+J_3+J_4.
\endaligned
\lab{(3.9)}
\ee

{\it Step 3.} It is clear that
$
|J_4| \le \frac{C}{r^3} \int_{B_r(x)-B_{r/2}(x)} |u(y, t)| dy.
$
The estimates $J_1, J_2, J_3$ are done on p381 of \cite{Z:4} following \cite{O:1}. This
proves of part (a) of the lemma.

{\it Proof of part (b).}
 Given $(w, l) \in Q_{r/4}(x, t)$, (\ref{(3.9)}) still holds when $(x, t)$ is replaced by $(w, l)$.  It is clear that, for $(w, l) \in Q_{r/4}(x, t)$,
\[
|\nabla_w J_4|  \le \frac{C}{r^4} \int_{B_r(x)-B_{r/2}(x)} |u(y, t)| dy.
\] The rest of the proof of the lemma is the same as that on p382-285 of \cite{Z:4} by taking the gradient of $J_1, J_2,  J_3$. \qed

With the lemma in hand, the proof of Theorem \ref{th:1.2} is then the same as the scaled iteration on p385-388 of \cite{Z:4} with the additional mild integral terms $J_4$ and $|\nabla J_4|$. Note these additional terms are spatial integrals only. These account for the difference
between the result of Navier Stokes equation and related result for the heat equation.

{\bf Acknowledgement} We wish to sincerely thank Professor Hongjie Dong for pointing out  the issue of the theorem.

\bigskip

\noindent e-mail: qizhang@math.ucr.edu

\enddocument